\documentclass[11pt,twoside]{article}
\usepackage{latexsym}
\usepackage{psfrag}
\usepackage[hyperfootnotes=false]{hyperref}
\usepackage{amsfonts} 
\usepackage{amsmath} 
\usepackage{amssymb}
\usepackage{bbm}
\usepackage{amsthm} 
\usepackage{mathrsfs} 
\usepackage{comment}
\usepackage{url}
\usepackage{graphicx}
\usepackage[sort&compress]{natbib}
\usepackage{rotating}
\usepackage[config, labelfont={normalsize,bf}, textfont=normalsize]{caption,subfig}
\usepackage{epsfig}
\usepackage[show]{simple_notes}
\usepackage[stable]{footmisc}
\usepackage{endnotes} 
\usepackage[english]{babel}

\theoremstyle{definition}

\theoremstyle{plain}
\newtheorem{theorem}{Theorem}
\newtheorem{proposition}{Proposition}

\newtheorem{lemma}{Lemma}

\theoremstyle{definition}

\newtheorem{example}{Example}

\theoremstyle{remark}

\newcommand{\R}{\mathbb{R}}

\renewcommand{\hat}[1]{\widehat{#1}}

\newcommand{\tsum}{\textstyle\sum}

\renewcommand{\P}{\mathbb{P}}
\newcommand{\E}{\mathbb{E}}

\newcommand{\e}{\epsilon}

\newcommand{\tr}{\operatorname{tr}}
\newcommand{\var}{\operatorname{Var}}

\newcommand{\hh}{{\bf{H}}}
\newcommand{\mnorm}[1]{\left\vert\kern-1.5pt\left\vert\kern-1.5pt\left\vert #1\right\vert\kern-1.5pt\right\vert\kern-1.5pt\right\vert}

    \def\independenT#1#2{\mathrel{\setbox0\hbox{$#1#2$}%
    \copy0\kern-\wd0\mkern4mu\box0}} 
\def\blfootnote{\xdef\@thefnmark{}\@footnotetext}
\long\def\symbolfootnote[#1]#2{\begingroup%
\def\thefootnote{\fnsymbol{footnote}}\footnote[#1]{#2}\endgroup} 

\usepackage{fullpage}
\newcommand{\mydef}{\ensuremath{: \, =}}

\makeatletter
\long\def\@makecaption#1#2{
        \vskip 0.8ex
        \setbox\@tempboxa\hbox{\small {\bf #1:} #2}
        \parindent 1.5em  
        \dimen0=\hsize
        \advance\dimen0 by -3em
        \ifdim \wd\@tempboxa >\dimen0
                \hbox to \hsize{
                        \parindent 0em
                        \hfil 
                        \parbox{\dimen0}{\def\baselinestretch{0.96}\small
                                {\bf #1.} #2
                                } 
                        \hfil}
        \else \hbox to \hsize{\hfil \box\@tempboxa \hfil}
        \fi
        }

\begin{document}

\begin{center}

{\bf{\LARGE{A More Powerful Two-Sample Test in High Dimensions using
Random Projection}}}

\vspace*{.3in}

{\large{
\begin{tabular}{ccc}
Miles E. Lopes$^{1}$ & Laurent J. Jacob$^1$ &  Martin J. Wainwright$^{1,2}$\\
{\footnotesize{{\texttt{mlopes@stat.berkeley.edu}}}} & {\footnotesize{{\texttt{laurent@stat.berkeley.edu}}}}  & {\footnotesize{{\texttt{wainwrig@stat.berkeley.edu}}}} 
\end{tabular}


\vspace*{.2in}
Departments of Statistics$^1$ and EECS$^2$\\
University of California, Berkeley
%
}}

\vspace*{.1in}

\begin{abstract}

 \noindent We consider the hypothesis testing problem of detecting a
 shift between the means of two multivariate normal distributions in
 the high-dimensional setting, allowing for the data dimension $p$ to
 exceed the sample size $n$. Specifically, we propose a new test
 statistic for the two-sample test of means that integrates a random
 projection with the classical Hotelling $T^2$ statistic.  Working
 under a high-dimensional framework with $(p,n)\to\infty$, we first
 derive an asymptotic power function for our test, and then provide
 sufficient conditions for it to achieve greater power than other
 state-of-the-art tests. Using ROC curves generated from synthetic
 data, we demonstrate superior performance against competing
 tests in the parameter regimes anticipated by our theoretical results.
 Lastly, we illustrate an advantage of our procedure's false positive
 rate with comparisons on high-dimensional gene expression data
 involving the discrimination of different types of cancer.

\end{abstract}

\end{center}

\section{Introduction}
Application domains such as molecular
biology and
fMRI~\cite[e.g.,][]{Lu2005Hotelling,Goeman2007Analyzing,fMRI1,fMRI2} 
have stimulated considerable interest in two-sample hypothesis testing problems in 
the high-dimensional setting, where two samples of data $\{X_1,\dots,X_{n_1}\}$ and
$\{Y_1,\dots,Y_{n_2}\}$ are subsets of $\R^p$, and $n_1,n_2 \ll p$.
The problem of discriminating between two data-generating
distributions becomes difficult in this context as the cumulative
effect of variance in many variables can ``explain away'' the correct
hypothesis. In transcriptomics, for instance, $p$ gene expression measures 
on the order of hundreds or thousands may be
used to investigate differences between two biological conditions, and
it is often difficult to obtain sample sizes $n_1$ and $n_2$ larger
than several dozen in each condition. For problems such as these,
classical methods may be ineffective, or not applicable at
all. Likewise, there has been growing interest in developing testing
procedures that are better suited to deal with the effects of
dimension \cite[e.g.,][]{BS96,SD2008,Srivastava2009,
  CQ2010,AUCopt}. 

A fundamental instance of the general two-sample
problem is the two-sample test of means with Gaussian data. In
this case, two independent sets of samples $\{X_1,\dots,X_{n_1}\}$ and
$\{Y_1,\dots,Y_{n_2}\}\subset \R^p$ are generated in an i.i.d. manner from
$p$-dimensional multivariate normal distributions $N(\mu_1,\Sigma)$
and $N(\mu_2,\Sigma)$ respectively, where the mean vectors
$\mu_1$ and $\mu_2$, and positive-definite covariance matrix $\Sigma\succ 0$, are all
fixed and unknown.
The hypothesis testing problem of interest is
\begin{equation}\label{origproblem}
{\bf{H}}_0: \mu_1=\mu_2 \text{ \ versus \ } {\bf{H}}_1: \mu_1\neq\mu_2.
\end{equation}
 The most well-known test statistic for this problem is the Hotelling
$T^2$ statistic, defined by
\begin{equation}
\label{Hotelling}
T^2:=\frac{n_1\, n_2}{n_1+n_2} \,
(\bar{X}-\bar{Y})^{\top}\hat{\Sigma}^{-1}\, (\bar{X}-\bar{Y}),
\end{equation}
where $\bar{X} := \frac{1}{n_1} \tsum_{j=1}^{n_1} X_j$ and
$\bar{Y} := \frac{1}{n_2} \tsum_{j=1}^{n_2} Y_j$ 
are the sample
means, and $\hat{\Sigma}$ is the pooled sample covariance matrix,
given by
{\small{$\hat{\Sigma} \mydef \textstyle{\frac{1}{n}\,
\tsum_{j=1}^{n_1}(X_j-\bar{X})(X_j-\bar{X})^{\top}+
\ \frac{1}{n}\tsum_{j=1}^{n_2}(Y_j-\bar{Y})(Y_j-\bar{Y})^{\top}}$}}, where
we define $n:=n_1+n_2-2$ for convenience.

When $p>n$, the matrix $\hat{\Sigma}$ is
singular, and the Hotelling test is not well-defined. Even when $p\leq n$,
the Hotelling test is known to perform poorly if $p$ is nearly as large as $n$.
This was shown in an important paper of Bai and Saranadasa (abbreviated BS)
 \cite{BS96}, who studied the performance 
of the Hotelling test under $(p,n) \to \infty$ with $p/n\to 1-\e$, and showed that 
the asymptotic power of the test suffers for small values of
$\e>0$. Consequently, several improvements on the Hotelling test have been proposed in the
high-dimensional setting in past years
\cite[e.g.,][]{BS96,SD2008,Srivastava2009,CQ2010}.

Due to the well-known degradation of $\hat{\Sigma}$ as an estimate of
$\Sigma$ in high dimensions, the line of research on extensions of
Hotelling test for problem \eqref{origproblem} has focused on
replacing $\hat{\Sigma}$ in the definition of $T^2$ with other
estimators of $\Sigma$. In the paper \cite{BS96}, BS proposed a test
statistic based on the quantity
$(\bar{X}-\bar{Y})^{\top}(\bar{X}-\bar{Y})$, which can be viewed as
replacing $\hat{\Sigma}$ with $I_{p\times p}$. It was shown by BS that
this statistic achieves non-trivial asymptotic power whenever the
ratio $p/n$ converges to a constant $c\in (0,\infty)$.
This statistic was later refined by Chen and Qin~\cite{CQ2010} (CQ
for short) who showed that the same asymptotic power can be achieved
without imposing any explicit restriction on the limit of
$p/n$. Another direction was considered by Srivastava and
Du~\cite{SD2008,Srivastava2009} (SD for short), who proposed
a test statistic based on $(\bar{X}-\bar{Y})^{\top}
\hat{D}^{-1}(\bar{X}-\bar{Y})$, where $\hat{D}$ is the diagonal matrix
associated with $\hat{\Sigma}$, i.e. $\hat{D}_{ii}=\hat{\Sigma}_{ii}$.
This choice ensures that $\hat{D}$ is invertible for all dimensions
$p$ with probability 1.
Srivastava and Du demonstrated that their test has superior asymptotic
power to the tests of BS and CQ under a certain parameter setting and
local alternative when $n=\mathcal{O}(p)$.  To the best of our
knowledge,
the procedures of CQ and SD represent the
state-of-the-art among tests for problem \eqref{origproblem}\footnote{The tests of BS, CQ, and SD 
actually extend somewhat beyond the problem \eqref{origproblem} in that their asymptotic power functions have been
obtained under data-generating distributions more general than Gaussian, e.g. satisfying simple moment conditions.}
with a known asymptotic power function
under the scaling $(p,n)\to\infty$.

In this paper, we propose a new testing procedure for problem
\eqref{origproblem} in the high-dimensional setting, which involves
randomly projecting the $p$-dimensional samples into a space of lower
dimension $k\leq \min\{n,p\}$, and then working with the Hotelling
test in $\R^k$.
Allowing $(p,n)\to\infty$, we derive an asymptotic power function for our 
test and show that it outperforms the tests of BS, CQ, and SD in 
terms of asymptotic relative efficiency under certain conditions. 
Our comparison results are valid with $p/n$ tending
to a constant or infinity. Furthermore, whereas the mentioned testing procedures can
only offer approximate level-$\alpha$ critical values,
our procedure specifies \emph{exact} level-$\alpha$ critical values
for general multivariate normal data. 
Our test is also very easy to implement, and has a computational cost
of order $\mathcal{O}(n^2p)$ operations when $k$ scales linearly with
$n$, which is modest in the high-dimensional setting.

From a conceptual point of view, the procedure studied here is most
distinct from past approaches in the way that covariance structure is
incorporated into the test statistic. As stated above, the test
statistics of BS, CQ, and SD are essentially based on versions of the
Hotelling $T^2$ with diagonal estimators of $\Sigma$. Our analysis and
simulations show that this limited estimation of $\Sigma$ sacrifices
power when the data variables are correlated, or when most of the
variance can be captured in a small number of variables.
In this regard, our procedure is motivated by the idea that
covariance structure may be used more effectively by testing with
projected samples in a space of lower dimension. The use of
projection-based test statistics has also been considered previously
in Jacob et al.  \cite{Jacob2010Gains} and Cl\'emen\c{c}on et
al. \cite{AUCopt}.

The remainder of this paper is organized as follows.  In
Section~\ref{SecBack}, we discuss the intuition for our testing
procedure, and then formally define the test
statistic. Section~\ref{SecMain} is devoted to a number of theoretical
results about the performance of the test. Theorem \ref{ThmPower} in
Section~\ref{SecPower} provides an asymptotic power function, and
Theorems \ref{ThmCQ} and \ref{ThmSD} in Sections \ref{sec: cqcomp} and
\ref{sec: sdcomp} give sufficient conditions for achieving greater
power than the tests of CQ and SD in the sense of asymptotic relative
efficiency.  In Sections~\ref{sec:ROC} and~\ref{sec:calibration}, we use synthetic data to make
performance comparisons with ROC and calibration curves against the mentioned tests,
as well as some recent non-parametric procedures such as maximum mean
discrepancy (MMD) \citep{Gretton2007A}, kernel Fisher discriminant
analysis (KFDA) \citep{Harchaoui2008Testing}, and a test based on
area-under-curve maximization, denoted TreeRank
\citep{AUCopt}. These simulations show that our test outperforms
competing tests in the parameter regimes anticipated by our theoretical
results. Lastly, in Section~\ref{sec:gene} we study an example
involving high-dimensional gene expression 
data, and demonstrate an advantage of our test in terms of its false positive
rate when discriminating between different types of cancer.
\paragraph{Notation.}We use $\delta:=\mu_1-\mu_2$ to denote the 
\emph{shift vector} between the distributions $N(\mu_1,\Sigma)$ and
$N(\mu_2,\Sigma)$. For a positive-definite covariance matrix $\Sigma$,
  let $D_{\sigma}$ be the diagonal matrix obtained by setting the off-diagonal
  entries of $\Sigma$ to 0, and also define the associated correlation
  matrix $R:=D_{\sigma}^{-1/2}\Sigma D_{\sigma}^{-1/2}$. Let
  $z_{1-\alpha}$ denote the $1-\alpha$ quantile of the standard normal
  distribution, and let $\Phi$ be its cumulative distribution
  function.  If $A$ is a matrix in $\R^{p\times p}$,
  let $\mnorm{A}_2$ denote its spectral norm (maximum singular value),
  and define the Frobenius norm $\mnorm{A}_F:=\sqrt{\sum_{i,j}
    A_{ij}^2}$. When all the eigenvalues of $A$ are real, we denote
  them by $\lambda_{\min}(A)=\lambda_p(A)\leq \cdots \leq
  \lambda_1(A)=\lambda_{\max}(A)$.  If $A$ is positive-definite, we
  write $A\succ 0$, and $A\succeq 0$ if $A$ is positive
  semidefinite. We use the notation $f(n)\lesssim g(n)$ if there is some
  absolute constant $c\in (0,\infty)$ such that the inequality
  $f(n)\leq c \, n$ holds for all large $n$. If both $f(n)\lesssim
  g(n)$ and $g(n)\lesssim f(n)$ hold, then we write $f(n)\asymp
  g(n)$. The notation $f(n)=o(g(n))$ means $f(n)/g(n)\to0$ as
  $n\to\infty$. For two random variables $X$ and $Y$, equality in
  distribution is written as $X\stackrel{d}{=} Y$.
%
%
%
%
%
%
%
%

\section{Random projection method}
\label{SecBack}

For the remainder of the paper, we retain the setup for the two-sample
test of means \eqref{origproblem} with Gaussian data given in Section
1.  In particular, our procedure can be implemented with $p>n$ or
$p\leq n$, as long as $k$ is chosen such that $k\leq \min\{n,p\}$. In
Section \ref{SecProjDim}, we demonstrate an optimality property of the
choice $k=\lfloor n/2 \rfloor$, which is valid in moderate or
high-dimensions, i.e., $p\geq \lfloor n/2\rfloor$, and we restrict our
attention to this case in Theorems \ref{ThmCQ} and \ref{ThmSD}.
%
\subsection{Intuition for random projection method}
\label{sec:intuition}
At a high level, our method can be viewed as a two-step
procedure. First, a single random projection $P_k^{\top}\in
\R^{k\times p}$ is drawn, and is used to map the samples from the
high-dimensional space $\R^p$ to a low-dimensional space $\R^k$.
Second, the Hotelling $T^2$ test is applied to a new hypothesis
testing problem, denoted ${\bf{H}}_{0,\text{proj}}$ versus
${\bf{H}}_{1,\text{proj}}$, in the projected space. A decision is then
pulled back to the high-dimensional problem \eqref{origproblem} by
simply rejecting the original null hypothesis ${\bf{H}}_0$ whenever
the Hotelling test rejects ${\bf{H}}_{0,\text{proj}}$ in the projected
space.

To provide some intuition for our method, it is possible to consider
the problem \eqref{origproblem} in terms of a competition between the
dimension $p$, and the ``statistical distance'' separating $\hh_0$ and
$\hh_1$.  On one hand, the accumulation of variance from a large
number of variables makes it difficult to discriminate between the
hypotheses, and thus, it is desirable to reduce the dimension of the
data. On the other hand, methods for reducing dimension also tend to
bring $\hh_0$ and $\hh_1$ ``closer together," making them harder to
distinguish. Mindful of the fact that the Hotelling $T^2$ measures the
separation of $\hh_0$ and $\hh_1$ in terms of the Kullback-Leibler
divergence $D_{\text{KL}}(N(\mu_1,\Sigma)\| N(\mu_2,\Sigma)) =
\frac{1}{2}\delta^{\top}\Sigma^{-1}\delta$, with
$\delta=\mu_1-\mu_2$,\footnote{When $p\leq n$, the distribution of
  the Hotelling $T^2$ under both $\hh_0$ and $\hh_1$ is given by a
  scaled noncentral $F$ distribution $\frac{p\,
    n}{n-p+1}\,F_{p,n-p-1}(\eta)$, with noncentrality parameter
  $\eta:=\frac{n_1\,n_2}{n_1+n_2}\, \delta^{\top}\Sigma^{-1}\delta$.
  The expected value of $T^2$ grows linearly with $\eta$,
  e.g., see Muirhead \cite[p. 216, p. 25]{Muirhead}.}  we see
that the relevant statistical distance is driven by the
length of $\delta$. Consequently, we seek to transform the data in a
way that reduces dimension and preserves most of the length of
$\delta$ upon passing to the transformed distributions.  From this
geometric point of view, it is natural to exploit the fact that random
projections can simultaneously reduce dimension and approximately
preserve length with high probability \cite{randProj}.

In addition to reducing dimension in a way that tends to preserve statistical
distance between $\hh_0$ and $\hh_1$, random projections have two
other interesting properties with regard to the design of test
statistics. Note that when the Hotelling test statistic is constructed from the 
projected samples in a space of dimension
$k\leq \min\{n,p\}$, it is proportional to
$[P_k^{\top} (\bar{X}-\bar{Y})]^{\top}(P_k^{\top} \hat{\Sigma}
P_k)^{-1} [P_k^{\top}(\bar{X}-\bar{Y})]$.\footnote{For the choice of $P_k^{\top}$
given in Section \ref{formal}, the matrix $P_k^{\top} \hat{\Sigma}
P_k$ is invertible with probability 1.} Thus, whereas the tests of BS, CQ, and
SD replace $\hat{\Sigma}$ in the definition of $T^2$ with a diagonal
estimator, our procedure uses $P_k^{\top}\hat{\Sigma} P_k$ as a $k\times k$
surrogate for $\hat{\Sigma}$. The key advantage is that
$P_k^{\top}\hat{\Sigma} P_k$ retains some information about the
off-diagonal entries of $\Sigma$. Another benefit offered by random
projection concerns the robustness of critical values. In the
classical setting where $p\leq n$, the critical values of the
Hotelling test are exact in the presence of Gaussian data. It is also
well-known from the projection pursuit literature that the empirical
distribution of randomly projected data tends to be approximately
Gaussian \cite{FreedmanDiaconis}. Our procedure leverages these two
facts by first ``inducing Gaussianity'' and then applying a test that
has exact critical values for Gaussian data. Consequently, we expect
that the critical values of our procedure may be accurate even when
the $p$-dimensional data are not Gaussian, and this idea is
illustrated by a simulation with data generated from a mixture model,
as well as an example with real transcriptomic data in Section
\ref{sec:gene}. \\
\subsection{Formal testing procedure}\label{formal}
For an integer $k\in \{1,\dots,\min\{n,p\}\}$, let $P_k^{\top}\in
\R^{k\times p}$ denote a random matrix with i.i.d. $N(0,1)$
entries,\footnote{We refer to $P_k^{\top}$ as a projection, even
  though it is not a projection in the strict sense of being
  idempotent. Also, we do not normalize $P_k^{\top}$ by $1/\sqrt{k}$
  (which is commonly used for Gaussian matrices \cite{randProj})
  because our statistic $T_k^2$ is invariant with respect to this
  scaling.} drawn independently of the data. Conditioning on a given
draw of $P_k^{\top}$, the projected samples $\{P_k^{\top} X_1,\dots,
P_k^{\top} X_{n_1}\}$ and $\{P_k^{\top} Y_1,\dots, P_k^{\top}
Y_{n_2}\}$ are distributed i.i.d. according to $N(P_k^{\top}\mu_i,
P_k^{\top}\Sigma P_k)$ respectively, with $i=1,2$. Since the projected
data are Gaussian and lie in a space of dimension no larger than $n$,
it is natural to consider applying the Hotelling test to the following
two-sample problem in the projected space $\R^k$:
\begin{equation}\label{projproblem}
{\bf{H}}_{0,\text{proj}}: P_k^{\top}\mu_1=P_k^{\top}\mu_2 \text{ \ versus \ } {\bf{H}}_{1,\text{proj}}: P_k^{\top} \mu_1\neq P_k^{\top}\mu_2.
\end{equation} 
For this projected problem, the Hotelling test statistic takes the 
form
\[
T_k^2:= \frac{n_1n_2}{n_1+n_2}[P_k^{\top} (\bar{X}-\bar{Y})]^{\top}(P_k^{\top}
\hat{\Sigma} P_k)^{-1} [P_k^{\top}(\bar{X}-\bar{Y})], 
\]
where $\bar{X}$, $\bar{Y}$, and $\hat{\Sigma}$ are as stated in the introduction.
Note that $P_k^{\top}\hat{\Sigma} P_k$ is invertible with probability 
1 when $P_k^{\top}$ has i.i.d. $N(0,1)$ entries, which ensures
that $T_k^2$ is well-defined, even when $p>n$.

When conditioned on a draw of $P_k^{\top}$, the $T_k^2$ statistic has an
$\frac{k\, n}{n-k+1}F_{k,n-k+1}$ distribution under
$\hh_{0,\text{proj}}$, since it is an instance of the Hotelling test statistic \cite[p. 216]{Muirhead}.
Inspection of the formula for $T_k^2$ also shows that its distribution
is the same under both ${\bf{H}}_0$ and
${\bf{H}}_{0,\text{proj}}$. Consequently, if we let
$t_{\alpha}:=\frac{k \, n}{n-k+1} F^{1-\alpha}_{k,n-k+1}$, where
$F^{1-\alpha}_{k,n-k+1}$ is the $1-\alpha$ quantile of the
$F_{k,n-k+1}$ distribution, then the condition $T_k^2 \geq t_{\alpha}$
is a level-$\alpha$ decision rule for rejecting the null hypothesis in
both the projected problem~\eqref{projproblem} and the original
problem~\eqref{origproblem}. Accordingly, we \emph{define} this as the
condition for rejecting ${\bf{H}}_0$ at level $\alpha$ in our
procedure for \eqref{origproblem}.  We summarize the implementation
of our procedure below.
%
\begin{center}
Implementation of random projection-based test at level
  $\alpha$ for problem \eqref{origproblem}.
\end{center}
\begin{equation*}\tag{$\star$}\label{procedure}
\boxed{
\begin{split}
&\text{1. Generate a single random matrix $P_k^{\top} \in \R^{k\times p}$
    with i.i.d. $N(0,1)$ entries.}\\ &\text{2. Compute $T_k^2$, using
    $P_k^{\top}$ and the two sets of samples.}\\ &\text{3. If
    $T_k^2\geq t_{\alpha}$, reject ${\bf{H}}_0$; otherwise accept
    ${\bf{H}}_0$.}\\
\end{split}
}
\end{equation*}

%
%
%
%
%
%
%
%

\section{Main results and their consequences}
\label{SecMain}

This section is devoted to the statement and discussion of our main
theoretical results, including an asymptotic power function for our test 
(Theorem~\ref{ThmPower}), and comparisons of asymptotic relative 
efficiency with state-of-the-art tests proposed in past work 
(Theorems~\ref{ThmCQ} and~\ref{ThmSD}).

\subsection{Asymptotic power function}
\label{SecPower}
%
%
%
%
%
%
%
%
%
%
%
Our first main result characterizes the asymptotic power of the
$T_k^2$ test statistic in the high-dimensional setting. As is standard
in high-dimensional asymptotics, we consider a sequence of hypothesis
testing problems indexed by $n$, allowing the dimension $p$, sample
sizes $n_1$ and $n_2$, mean vectors $\mu_1$ and $\mu_2$ and covariance
matrix $\Sigma$ to implicitly vary as functions of $n$, with $n$
tending to infinity. We also make another type of asymptotic
assumption, known as a \emph{local alternative}, which is commonplace
in hypothesis testing (e.g., see van der Vaart~\cite[\S14.1]{vdv}).
The idea lying behind a local alternative assumption is that if the
difficulty of discriminating between $\hh_0$ and $\hh_1$ is ``held
fixed'' with respect to $n$, then it is often the case that most
testing procedures have power tending to 1 under $\hh_1$ as
$n\to\infty$. In such a situation, it is not possible to tell if one
test has greater asymptotic power than another. Consequently, it is
standard to derive asymptotic power results under the extra condition
that $\hh_0$ and $\hh_1$ become harder to distinguish as $n$
grows. This theoretical device aids in identifying the conditions
under which one test is more powerful than another.  The following
local alternative ({\bf{A0}}), and balancing assumption ({\bf{A1}}),
are the same as those used by Bai and Saranadasa~\cite{BS96} to study
the asymptotic power of the classical Hotelling test under
$(n,p)\to\infty$.
In particular, the local alternative ({\bf{A0}}) means that the
Kullback-Leibler divergence between the $p$-dimensional sampling
distributions, $D_{\text{KL}}(N(\mu_1,\Sigma)\, \| \,
N(\mu_2,\Sigma))=\frac{1}{2}\delta^{\top}\Sigma^{-1}\delta$, tends to
$0$ as $n\to\infty$.\\

({\bf{A0}}) (Local alternative.) The shift vector and covariance matrix satisfy $\delta^{\top}\Sigma^{-1}\delta =o(1)$.\\

({\bf{A1}}) There is a constant $b\in(0,1)$ such that $n_1/n \to b$.\\

({\bf{A2}}) There is a constant $y\in (0,1)$ such that $k/n\to y$.\\

To set some notation for our asymptotic power result in
Theorem~\ref{ThmPower}, let $\theta:=(\delta,\Sigma)$ be an ordered pair
containing the relevant parameters for problem \eqref{origproblem},  and define $\Delta_k^2$ as twice the
Kullback-Leibler divergence between the projected sampling
distributions,
\begin{align}
\label{KL}
\Delta_k^2 & := 2 \, D_{\text{KL}}\left(N(P_k^{\top}\mu_1,
P_k^{\top}\Sigma P_k) \, \big\| \, N(P_k^{\top}\mu_2, P_k^{\top}\Sigma
P_k)\right) \; = \; \delta^{\top}P_k(P_k^{\top}\Sigma
P_k)^{-1}P_k^{\top}\delta.
\end{align}
When interpreting the statement of Theorem \ref{ThmPower} below,
it is important to notice that each time the procedure
(\ref{procedure}) is implemented, a draw of $P_k^{\top}$ induces a new
test statistic $T_k^2$. Making this dependence on $P_k^{\top}$
explicit, let $\beta(\theta; P_k^{\top})$ denote the exact
(non-asymptotic) power function of the $T_k^2$ statistic at level
$\alpha$ for problem \eqref{origproblem}, conditioned on a given draw
of $P_k^{\top}$, as in procedure ($\star$).

\begin{theorem}
\label{ThmPower}
Assume conditions ({\bf{A0}}), ({\bf{A1}}), and ({\bf{A2}}). Then, for almost all sequences of projections $P_k^{\top}$, the
power function $\beta(\theta ;P_k^{\top})$ satisfies
\begin{align}
\label{asympower}
\beta(\theta ; P_k^{\top}) - \Phi\left(-z_{1-\alpha}+
b(1-b)\sqrt{\frac{1-y}{2y}} \, \Delta_{k}^2\,\sqrt{n} \, \right) & \to
0 \quad \mbox{ as $n \to \infty$.}
\end{align}
\end{theorem}
\noindent {\bf{Remarks.}}
Notice that if $\Delta_k^2=0$ (e.g. under $\hh_0$), then $\Phi(-z_{1-\alpha}+0)= \alpha$, which
corresponds to blind guessing at level $\alpha$. Consequently, the
second term $b(1-b)\sqrt{\frac{1-y}{2y}} \,
\Delta_{k}^2  \, \sqrt{n}$
determines the advantage of our procedure over blind guessing. Since
$\Delta_k^2$ is twice the KL-divergence between the
projected sampling distributions, these observations conform to the
intuition from Section 2 that the KL-divergence measures the
discrepancy between $\hh_0$ and $\hh_1$.\\

\noindent \emph{Proof of Theorem \ref{ThmPower}.} Let
$\beta_{\text{H}}(\theta; P_k^{\top})$ denote the exact power of the
Hotelling test for the projected problem \eqref{projproblem} at level
$\alpha$.  As a preliminary step, we verify that
\begin{equation}\label{equiv}
\beta(\theta;P_k^{\top})=\beta_{\text{H}}(\theta; P_k^{\top}),
\end{equation}
for almost all $P_k^{\top}$. To see this, first recall from Section 2 that 
the condition $T_k^2\geq t_{\alpha}$ is a level-$\alpha$ rejection
criterion in both the procedure \eqref{procedure} for the original
problem \eqref{origproblem}, and the Hotelling test for the projected
problem \eqref{projproblem}. Next, note that if ${\bf{H}}_1:
\delta\neq 0$ holds, then ${\bf{H}}_{1,\text{proj}}: P_k^{\top}\delta
\neq 0$ holds with probability 1, since $P_k^{\top}\delta$ is
distributed as $N(0,\|\delta\|_2^2 \, I_{k\times k})$. Consequently,
for almost all $P_k^{\top}$, the level-$\alpha$ decision rule
$T_k^2\geq t_{\alpha}$ has the same power against the alternative in
both the original and the projected problems, which verifies
\eqref{equiv}. This establishes a technical link that allows results
on the power of the classical Hotelling test to be transferred to the
high-dimensional problem \eqref{origproblem}.

In order to complete the proof, we use a result of Bai and Saranadasa
\cite[Theorem 2.1]{BS96}\footnote{To prevent confusion, note that the
  notation in BS \cite{BS96} for $\delta$ differs from ours.},
which asserts that if $\Delta_k^2 =o(1)$ holds for a fixed sequence of
projections $P_k^{\top}$, and assumptions ({\bf{A1}}) and ({\bf{A2}}) hold, then
$\beta_{\text{H}}(\theta;P_k^{\top})$ satisifes
\begin{equation}\label{KLlimit}
\beta_{\text{H}}(\theta; P_k^{\top}) - \Phi\left(-z_{1-\alpha}+
b(1-b)\sqrt{\frac{1-y}{2y}} \, \Delta_{k}^2\,\sqrt{n} \, \right)\to 0
\text{ \ as \ } n\to\infty.
\end{equation}
To ensure $\Delta_k^2=o(1)$, we appeal to a deterministic matrix inequality 
that follows from the proof of Lemma 3 in Jacob et al. \cite{Jacob2010Gains}.
Namely, for any full rank matrix $M^{\top}\in \R^{k\times p}$, and any $\delta\in \R^p$,
\begin{equation*}
\delta^{\top}M(M^{\top}\Sigma M)^{-1}M^{\top}\delta \, \leq \, \delta^{\top}\Sigma^{-1}\delta.
\end{equation*}
Since $P_k^{\top}$ is full rank with probability 1, we see that
$\Delta_k^2 \leq \delta^{\top}\Sigma^{-1}\delta \to 0$
for almost all sequences of $P_k^{\top}$ under the local alternative
({\bf{A0}}), as needed. Thus, the proof of Theorem \ref{ThmPower} is completed by
combining equation~\eqref{equiv} with the limit~\eqref{KLlimit}.\qed
%
%
%
%
%
%
%
%
\subsection{Asymptotic relative efficiency (ARE)}
\label{SecARE}

Having derived an asymptotic power function in Theorem~\ref{ThmPower},
we are now in position to provide a detailed comparison with the
tests of CQ~\cite{CQ2010} and SD~\cite{SD2008, Srivastava2009}.
We denote the asymptotic power
function of our level-$\alpha$ random projection-based test (RP) by
\begin{equation}\label{ProjPower}
\beta_{\text{RP}}(\theta;P_k^{\top}):=
\Phi\left(-z_{1-\alpha}+
b(1-b)\sqrt{\frac{1-y}{2y}} \, \Delta_{k}^2\,\sqrt{n} \,
\right),
\end{equation}
where we recall $\theta := (\delta,\Sigma)$. The asymptotic power functions for the level-$\alpha$ testing procedures of CQ~\cite{CQ2010} and SD
\cite{SD2008, Srivastava2009} are given by
\begin{subequations}\label{CQSDpower}
\label{EqnOtherPow}
\begin{align}
    \beta_{\text{CQ}}(\theta) & := \ \Phi\left(-z_{1-\alpha}+ {\frac{b
        (1-b) }{\sqrt{2}}\, \frac{\|\delta\|_2^2 \, n}{\mnorm{\Sigma}_F}}
    \right), \quad \mbox{and} \\
 \beta_{\text{SD}}(\theta)&:=
    \ \Phi\left(-z_{1-\alpha}+ {\frac{b(1-b)}{\sqrt{2}} \, \frac{
        \delta^{\top} D_{\sigma}^{-1}\delta \, n}{\mnorm{R}_F} }\right),
\end{align}
\end{subequations}
where $D_\sigma$ denotes the matrix formed by setting the
off-diagonal entries of $\Sigma$ to 0, and $R$ denotes the correlation
matrix associated to $\Sigma$.  The functions $\beta_{\text{CQ}}$ and
$\beta_{\text{SD}}$ are derived under local alternatives and
asymptotic assumptions that are similar to the ones used here to
obtain $\beta_{\text{RP}}$. In particular, all three functions can be
obtained allowing $p/n$ to tend to an arbitrary positive constant, or
to infinity.
%
%
%
%
%
%


  A standard method of comparing asymptotic power
functions is through the concept of
\emph{asymptotic relative efficiency}, or ARE for short (e.g.,
see van der Vaart~\cite[ch. 14-15]{vdv}). Since the term added to
$-z_{1-\alpha}$ inside the $\Phi$ function is what controls power, the
relative efficiency of tests is defined by the ratio of such
terms.  More explicitly, we define
\begin{subequations}
\begin{align}
\text{ARE}\left(\beta_{\text{CQ}}; \beta_{\text{RP}}\right)& : =
\Big(\textstyle{{\frac{\|\delta\|_2^2 \, n}{\mnorm{\Sigma}_F}}\Big /
\sqrt{\frac{1-y}{y}} \,\Delta_k^2\,\sqrt{n} }\Big)^2, \quad
\mbox{and} \\
\text{ARE}\left(\beta_{\text{SD}}; \beta_{\text{RP}}\right)&: =
\Big(\textstyle{{\frac{\delta^{\top}D_{\sigma}^{-1}\delta \, n}{\mnorm{R}_F}}
\Big/\sqrt{\frac{1-y}{y}}\,\Delta_k^2\,\sqrt{n} }\Big)^2.
\end{align}
\end{subequations}
Whenever the ARE is less than 1, our procedure is considered to have
greater asymptotic power than the competing test---with our advantage
being greater for smaller values of the ARE. Consequently, we seek
sufficient conditions in Theorems~\ref{ThmCQ} and \ref{ThmSD} below
for ensuring that the ARE is small.

In classical analyses of asymptotic relative efficiency, the ARE is usually
a deterministic quantity that 
does not depend on $n$. However, in the current context, our use of high-dimensional
asymptotics, as well as a randomly constructed
test statistic, lead to an ARE that varies with $n$ and is random.
(In other words, the ARE specifies a sequence of random variables indexed 
by $n$.)
Moreover, the dependence of the ARE on $\Delta_k^2$ implies that the
ARE is affected by the orientation of the shift vector $\delta$.\footnote{In
  fact, $\text{ARE}\left(\beta_{\text{CQ}}; \beta_{\text{RP}}\right)$
  and $\text{ARE}\left(\beta_{\text{SD}}; \beta_{\text{RP}}\right)$
  are invariant with respect to scaling of $\delta$, and so the
  orientation $\delta/\|\delta\|_2$ is the only part of the shift
  vector that is relevant for comparing power.} To consider an
average-case scenario, where no single orientation of $\delta$ is of
particular importance, we place a prior on $\delta$, and assume that
it follows a spherical distribution\footnote{
  i.e. $\delta\stackrel{d}{=}U\delta$ for any orthogonal matrix $U$.}
with $\P(\delta=0)=0$.  This implies that the orientation
$\delta/\|\delta\|_2$ of the shift follows the uniform (Haar)
distribution on the unit sphere. We
emphasize that our procedure (\ref{procedure}) does not rely on this
choice of prior, and that it is only a device for making an
average-case comparison against CQ and SD in Theorems~\ref{ThmCQ}
and~\ref{ThmSD}. Lastly, we point out that a similar assumption was
considered by Srivastava and Du~\cite{SD2008}, who let $\delta$ be a
deterministic vector with all coordinates equal to the same value, in
order to compare with the results of BS~\cite{BS96}.

To be clear about the meaning of Proposition \ref{PropKeyLimit} and
Theorems~\ref{ThmCQ} and \ref{ThmSD} below, we henceforth regard the
ARE as a function of two random objects, $P_k^{\top}$ and $\delta$,
and our probability statements are made with this understanding. We
complete the preparation for our comparison theorems by stating
Proposition~\ref{PropKeyLimit} and several limiting assumptions with
$n\to\infty$.\\

({\bf{A3}}) The shift $\delta$ has a spherical distribution with
$\P(\delta = 0)=0$, and is independent of $P_k^{\top}$.\\

({\bf{A4}}) There is a constant $a\in[0,1)$ such that $k/p \to a$.\\

({\bf{A5}}) Assume $\frac{1}{\sqrt{k}} \, \frac{\tr(\Sigma)}{p \,
    \lambda_{\text{min}}(\Sigma)} = o(1)$.\\

({\bf{A6}}) Assume
  $\frac{\mnorm{D_{\sigma}^{-1}}_2}{\tr(D_{\sigma}^{-1})}=o(1)$.\\
  
As can be seen from the formulas for $\beta_{\text{RP}}$ and the ARE,
the performance of the $T_k^2$ statistic is determined by the
random quantity $\Delta_k^2$. The following proposition provides
interpretable upper and lower bounds on $\Delta_k^2$ that hold with
high-probability. This proposition is the main technical tool needed
for our comparison results in Theorems~\ref{ThmCQ} and~\ref{ThmSD}. A
proof is given in Appendix~\ref{AppPropKeyLimit}.
%
%
%
\begin{proposition}
\label{PropKeyLimit}
Under conditions ({\bf{A3}}), ({\bf{A4}}), and ({\bf{A5}}), let $c$ be
any positive constant strictly less than $(1-\sqrt{a})^2$, and let $C$
be any constant strictly greater than
$\frac{(1+\sqrt{a})^2}{(1-\sqrt{a})^2}$. Then, as $n\to\infty$, we
have
\begin{subequations}
\begin{align}
\label{lowerProp}
\P \left( \frac{\Delta_k^2}{\|\delta\|_2^2} \geq \frac{c \,
    k}{\tr(\Sigma)} \right) & \to 1, \quad \mbox{and} \\
\bigskip
\label{upperProp}
\P \left( \frac{\Delta_k^2}{\|\delta\|_2^2} \leq \frac{C \, k}{p \,
    \lambda_{\min}(\Sigma)} \right) & \to 1.
\end{align}
\end{subequations}
\end{proposition}
\noindent {\bf{Remarks.}}  Although we have presented upper and lower
bounds in an asymptotic manner, our proof specifies non-asymptotic
bounds on $\Delta_k^2/\|\delta\|_2^2$. Due to the fact that
Proposition~\ref{PropKeyLimit} is a tool for making asymptotic
comparisons of power in Theorems \ref{ThmCQ} and \ref{ThmSD}, it is
sufficient and simpler to state the bounds in this asymptotic
form. Note that if the condition $\tr(\Sigma)\asymp p \,\lambda_{\min}(\Sigma)$ 
holds, then Proposition 1 is sharp in the
sense the upper and lower bounds \eqref{lowerProp} and
\eqref{upperProp} match up to constants.
%
%
%
%
\subsection{Choice of projection dimension $k=\lfloor n/2 \rfloor$ }
\label{SecProjDim}
We now demonstrate an optimality property of the choice of projected
dimension $k=\lfloor n/2\rfloor$. Note that this choice implicitly assumes 
$p\geq \lfloor n/2\rfloor$, but this does not affect the applicability of 
procedure \eqref{procedure} in moderate or high-dimensions.
Letting $k/n\to y\in (0,1)$ as in assumption ({\bf{A2}}), recall that the 
asymptotic power function from Theorem \ref{ThmPower} is
\begin{equation*}
\Phi\left(-z_{1-\alpha}+ b(1-b)\sqrt{\frac{1-y}{2y}} \,\Delta_{k}^2\,
\sqrt{n} \, \right).
\end{equation*}
Since Proposition~\ref{PropKeyLimit} indicates that $\Delta_k^2$
scales linearly in $k$ up to random fluctuations, we see that formally
replacing $k$ with $y \, n$ leads to maximizing the function $f(y):=
\sqrt{\frac{1-y}{2y}}\, y$. The fact that $f$ is maximized at $y=1/2$
suggests that in certain cases, $k=\lfloor n/2 \rfloor$ may be
asymptotically optimal in a suitable sense. Considering a simple case
where $\Sigma = \sigma^2 \, I_{p\times p}$ for some absolute constant
$\sigma^2>0$, it can be shown\footnote{Note that $\|\delta\|_2$ and $\delta/\|\delta\|_2$ are independent,
and $\E\left[\frac{\delta^{\top}}{\|\delta\|_2}A \frac{\delta}{\|\delta\|_2}\right] = \tr(A)/p$ for any $A\in \R^{p\times p}$, under ({\bf{A3}});
see \cite[p. 38]{Muirhead}.}
that under assumptions ({\bf{A2}}),
({\bf{A3}}), and integrability of $\|\delta\|_2^2$,
\begin{equation}\label{limitassumption}
\frac{p}{n}\,\frac{\E(\Delta_k^2)}{\E(\|\delta\|_2^2)}\to y/\sigma^2,
\end{equation}
for all \,$y\in (0,1)$, as $n\to\infty$.  The following proposition is
an immediate extension of this observation, and shows that $k=\lfloor
n/2\rfloor$ is optimal in a precise sense for parameter settings that
include $\Sigma=\sigma^2\, I_{p\times p}$ as a special case. Namely,
as $n\to\infty$, the quantity $\sqrt{\frac{1-y}{2y}}\,\Delta_k^2$ is
largest on average for $k=\lfloor n/2\rfloor$ among all choices of
$k$, under the conditions stated below.
\begin{proposition}
\label{PropChoiceofk}
In addition to assumptions ({\bf{A2}}) and ({\bf{A3}}), suppose that
$\|\delta\|_2^2$ is integrable. Also assume that for some absolute constant
$\sigma^2>0$, the limit~\eqref{limitassumption} holds for any $y\in
(0,1)$. Let \,$y^*=1/2$, and $k^*=\lfloor n/2 \rfloor$. Then, for any
\,$y\in (0,1)$,
\begin{equation}
\lim_{n\to\infty}\frac{\sqrt{\frac{1-y^*}{2y^*}}\,\E(\Delta_{k^*}^2)}{\sqrt{\frac{1-y}{2y}}\,\E(\Delta_{k}^2)}
= \frac{1}{2\sqrt{y(1-y)}} \geq 1.
\end{equation}
\begin{figure}[h]
\centering
\subfloat[Setting (1)  ]{\includegraphics[angle=270,
  width=.45\linewidth]{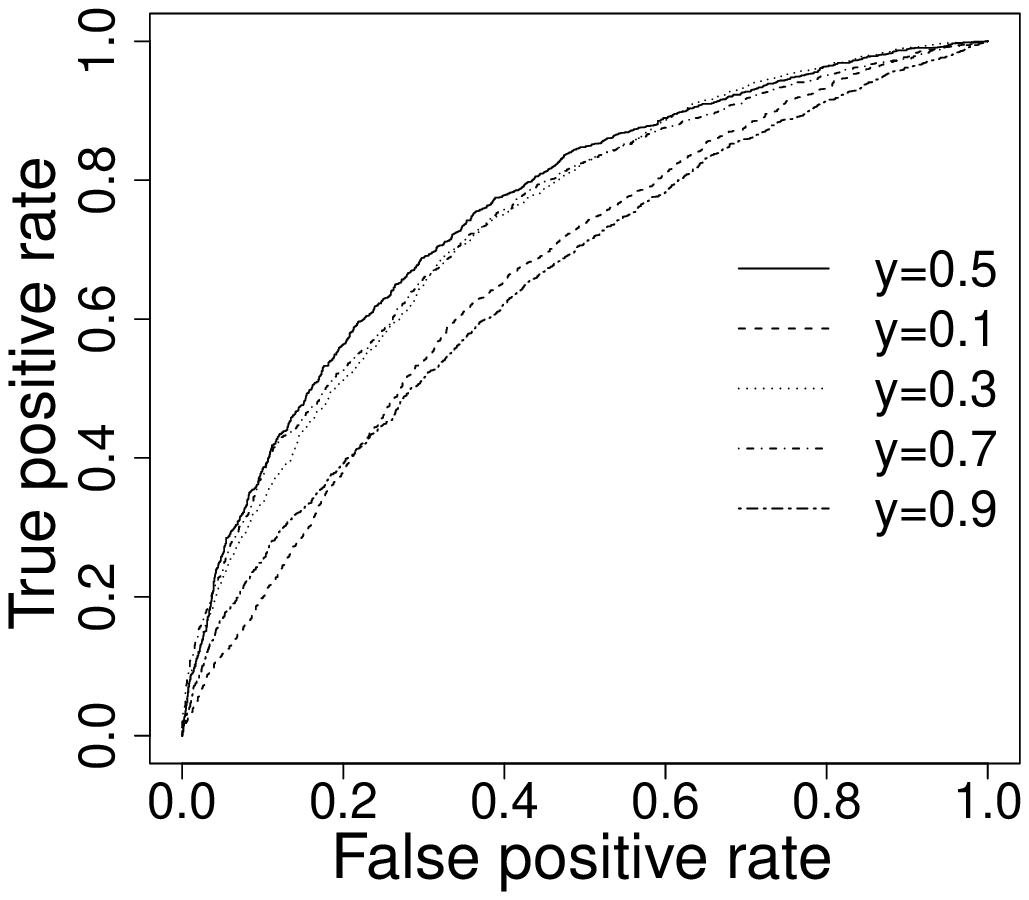}} \ \ \ \ 
\subfloat[Setting (2)  ]{\includegraphics[angle=270,
  width=.45\linewidth]{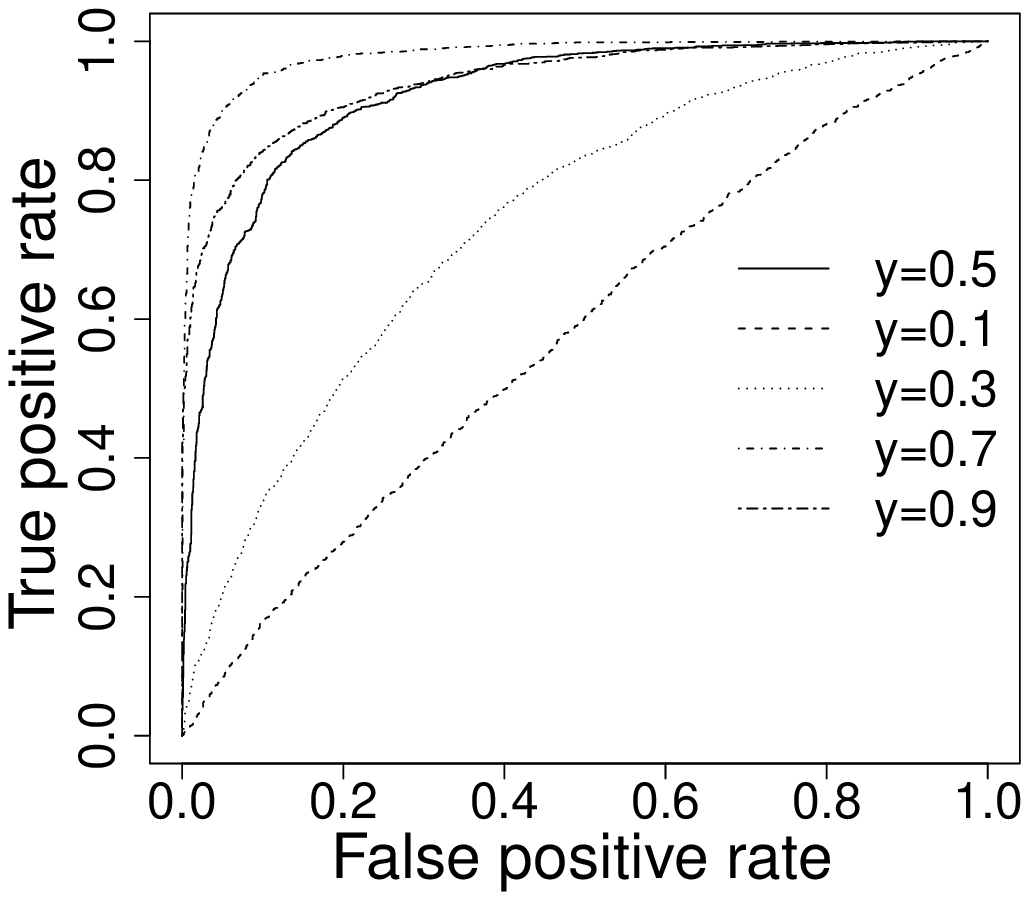}}
  \caption{Setting (1) corresponds to $\Sigma = \sigma^2 \, I_{p\times p}$ with $\sigma^2=50$, and Setting (2) involves a covariance matrix $\Sigma$ with randomly selected eigenvectors and a rapidly decaying spectrum.
  The ROC curves indicate that $k=\lfloor n/2\rfloor$ is optimal, or nearly optimal, among the five choices of $y$ in the two settings.}
\label{fig:sweep}
\end{figure}
\end{proposition}
\noindent {\bf{Remarks.}} The ROC curves in Figure \ref{fig:sweep}
illustrate several choices of projection dimension, with $k=\lfloor y
n \rfloor$ and $y=0.1, \, 0.3,\, 0.5,\,0.7,\,0.9$, under two different
parameter settings. In Setting (1), $\Sigma =\sigma^2\, I_{p\times p}$
with $\sigma^2=50$, and in Setting (2), the matrix $\Sigma$ was
constructed with a rapidly decaying spectrum, and a matrix of eigenvectors
drawn from the uniform (Haar) distribution on the orthogonal group,
as in panel (d) of Figure \ref{fig:exps} (see Section \ref{sec:ROC}
for additional details). The curves in both settings were generated by
sampling $n_1=n_2=50$ data at points from each of the distributions
$N(\mu_1,\Sigma)$ and $N(\mu_2,\Sigma)$ in $p=200$ dimensions, and
repeating the process 2000 times under both $\hh_0$ and $\hh_1$. For
the experiments under $\hh_1$, the shift $\delta$ was drawn uniformly
from a sphere of radius 3 for Setting (1), and radius 1 for Setting (2)---in 
accordance with assumption ({\bf{A3}}) in
Proposition~\ref{PropChoiceofk}. Note that $k=\lfloor n/2\rfloor$
gives the best ROC curve for  Setting (1) in Figure \ref{fig:sweep}, which agrees with the fact
that $\Sigma=\sigma^2 \, I_{p\times p}$ satisfies the conditions of
Proposition \ref{PropChoiceofk}. In Setting (2), we see that the
choice $k=\lfloor n/2\rfloor$ is not far from optimal, even when
$\Sigma$ is very different from $\sigma^2 I_{p\times p}$.

%
%
%
%
%
\subsection{Power comparison with CQ}
\label{sec: cqcomp}
The next result provides a sufficient condition for the $T_k^2$
statistic to be asymptotically more powerful than the test of CQ. A
proof is given at the end of this section (\ref{sec: cqcomp}).
\begin{theorem}
\label{ThmCQ}
Under the conditions of Proposition~\ref{PropKeyLimit}, suppose
 that we use a projection dimension $k=\lfloor n/2\rfloor$, where we
 assume $p\geq \lfloor n/2\rfloor$. Fix a number $\e_1>0$, and
 let $c_1(\e_1)$ be any constant strictly greater than $\frac{4}{\e_1
   (1-\sqrt{a})^4}$. If the condition
\begin{equation}
\label{suffCQ}
n\geq c_1(\e_1) \, \frac{\tr(\Sigma)^2}{\mnorm{\Sigma}_F^2},
\end{equation}
holds for all large $n$, then \emph{$\P\left
  [{\text{ARE}}\left(\beta_{\text{CQ}}; \beta_{\text{RP}}\right) \leq
    \e_1 \right ] \to 1$} as $n\to\infty$.
\end{theorem}
\noindent {\bf{Remarks.}} The case of $\e_1=1$ serves as the reference
for equal asymptotic performance, with values $\e_1<1$ corresponding
to the $T_k^2$ statistic being asymptotically more powerful than the
test of CQ.  
To interpret the result, note that Jensen's inequality implies that
the ratio $\tr(\Sigma)^2/\mnorm{\Sigma}_F^2$ lies between $1$ and $p$,
for any choice of $\Sigma$. As such, it is reasonable to interpret
this ratio as a measure of the \emph{effective dimension} of the
covariance structure.\footnote{This ratio has also been studied as an
  effective measure of matrix rank in the context of low-rank matrix
  reconstruction~\cite{ellStar}.}
    The message of
Theorem~\ref{ThmCQ} is that as long as the sample size $n$ grows
faster than the effective dimension, then our projection-based test is
asymptotically superior to the test of CQ.
 
The ratio $\tr(\Sigma)^2/\mnorm{\Sigma}_F^2$ can also be viewed as
measuring the \emph{decay rate} of the spectrum of $\Sigma$, with the
condition $\tr(\Sigma)^2\big/\mnorm{\Sigma}_F^2 \ll p$ indicating
rapid decay.  This condition means that the data has low variance in
``most'' directions in $\R^p$, and so projecting onto a random set of
$k$ directions will likely map the data into a low-variance subspace
in which it is harder for chance variation to explain away the correct
hypothesis, thereby resulting in greater power.
 %
 %
 %

\begin{example}\label{spiked}
One instance of spectrum decay occurs when the top $s$ eigenvalues of
$\Sigma$ contain most of the mass in the spectrum.  When $\Sigma$ is
diagonal, this has the interpretation that $s$ variables capture most
of the total variance in the data. For simplicity, assume
$\lambda_1=\cdots =\lambda_s>1$ and $\lambda_{s+1}=\cdots
=\lambda_p=1$, which is similar to the \emph{spiked covariance model}
introduced by Johnstone \cite{Johnstone}.  If the top $s$ eigenvalues
contain half of the total mass of the spectrum, then $s \, \lambda_1 =
(p-s)$, and a simple calculation shows that
\begin{equation}
\begin{split}
\frac{\tr(\Sigma)^2}{\mnorm{\Sigma}_F^2} &= \frac{4 \, \lambda_1^2}{\lambda_1^2+\lambda_1}\, s \leq 4s.\\
\end{split}
\end{equation}
This again illustrates the idea that condition \eqref{suffCQ} is satisfied as 
long as $n$ grows at a faster rate than the effective number of variables $s$.
It is straightforward to check that this example satisfies assumption
({\bf{A5}}) of Theorem \ref{ThmCQ} when, for instance, $\lambda_1=o(\sqrt{k})$.
 \end{example}
 
 \begin{example}\label{decayEx}
Another example of spectrum decay can be specified by $\lambda_i(\Sigma) \propto  i^{-\nu}$, 
 for some absolute proportionality constant, a rate parameter $\nu \in (0,\infty)$, and $i=1,\dots, p$. 
 This type of decay arises in connection with the Fourier coefficients of functions in Sobolev ellipsoids
 \cite[\S 7.2]{Was06}. Noting that
 $\tr(\Sigma) \asymp \int_1^p x^{-\nu} dx$ 
 and
 $\mnorm{\Sigma}_F^2 \asymp \int_1^p x^{-2\nu}dx,$
direct computation of the integrals shows that
\begin{equation*}
\frac{\tr(\Sigma)^2}{\mnorm{\Sigma}_F^2}\asymp
\begin{cases} 1& \text{ if } \nu>1
\\
\log^2 p & \text{ if } \nu=1\\
p^{2(1-\nu)} & \text{ if } \nu\in (\frac{1}{2},1)\\
p/\log p & \text{ if } \nu=\frac{1}{2}\\
p & \text{ if } \nu \in (0,\frac{1}{2})\\
\end{cases}.
\end{equation*}
Thus, a decay rate given by $\nu\geq 1$ is easily sufficient for
condition \eqref{suffCQ} to hold unless the dimension grows
exponentially with $n$.  On the other hand, decay rates associated to
$\nu\leq 1/2$ are too slow for condition \eqref{suffCQ} to hold when
$n\ll p$, and rates corresponding to $\nu\in (\frac{1}{2},1)$ lead to
a more nuanced competition between $p$ and $n$. Assumption ({\bf{A5}})
of Theorem~\ref{ThmCQ} holds for all $\nu \in (0,1)$, but when $\nu=1$ or $\nu>1$, the
dimension $p$ must satisfy the extra conditions $\log p = o(\sqrt{k})$
or $p^{\nu-1} = o(\sqrt{k})$ respectively.\footnote{It may be possible to relax ({\bf{A5}}) with a more
  refined analysis of the proof of Proposition~\ref{PropKeyLimit}.} 
 \end{example}

\noindent The proof of Theorem~\ref{ThmCQ}\, is a direct application of
Proposition~\ref{PropKeyLimit}.\\
 
\noindent \emph{Proof of Theorem \ref{ThmCQ}.}  Recalling
$\text{ARE}\left(\beta_{\text{CQ}}; \beta_{\text{RP}}\right) =
\Big({\frac{n \|\delta\|_2^2}{\mnorm{\Sigma}_F}}\Big/\sqrt{n}
\Delta_k^2\Big)^2$, with $k=\lfloor n/2\rfloor$ and $y=1/2$, the event of interest,
\begin{equation}\label{first}
\text{ARE}\left(\beta_{\text{CQ}}; \beta_{\text{RP}}\right)\leq \e_1,
\end{equation}
is the same as
$$ \frac{n}{\mnorm{\Sigma}_F^2} \frac{1}{\e_1}\leq
\left(\frac{\Delta_k^2}{\|\delta\|_2^2}\right)^2.$$ By
Proposition~\ref{PropKeyLimit}, we know that for any positive constant
$c$ strictly less than $(1-\sqrt{a})^2$, the probability of the event
\begin{equation}\label{second}
\frac{c\, k}{\tr(\Sigma)}\leq \frac{\Delta_k^2}{\|\delta\|_2^2}.
\end{equation}
tends to 1 as $n\to\infty$. Consequently, as long as the inequality
\begin{equation} \label{determ}
 \frac{n}{\mnorm{\Sigma}_F^2}  \frac{1}{\e_1}\leq \left(\frac{c\, k}{\tr(\Sigma)}\right)^2,
\end{equation}
holds for all large $n$,
then the event \eqref{first} of interest will also have probability tending to 1. 
Replacing $k$ with $\frac{n}{2}\cdot [1-o(1)]$, the last condition is the same as
\begin{equation}
n\geq \frac{\tr(\Sigma)^2}{\mnorm{\Sigma}_F^2}\cdot \frac{4}{\e_1\, c^2\, [1-o(1)]^2}.
\end{equation}
Thus, for a given choice of $c_1(\e_1)$ in the statement of the theorem, 
it is possible to choose a positive $c<(1-\sqrt{a})^2$ so that inequality \eqref{determ} 
is implied by the claimed sufficient
condition \eqref{suffCQ} for all large $n$.\qed
%
%
%
%
%
%
%
%
%
\subsection{Power comparison with SD}
\label{sec: sdcomp}
We now give a sufficient condition for our procedure to be
asymptotically more powerful than SD.
\begin{theorem}\label{ThmSD}
In addition to the conditions of Theorem~\ref{ThmCQ}, assume that
({\bf{A6}}) holds.  Fix a number $\e_1>0$, and let $c_1(\e_1)$ be any
constant strictly greater than $\frac{4}{\e_1 (1-\sqrt{a})^4}$. If the
condition
\begin{equation}
\label{suffSD}
 n\geq c_1(\e_1) \left(\frac{\tr(\Sigma)}{p}\right)^2 
\left(\frac{\tr(D_{\sigma}^{-1})} { \mnorm{R}_F} \right)^2
\end{equation}
holds for all large $n$, then \emph{$\P\left
  [\text{ARE}\left(\beta_{\text{SD}}; \beta_{\text{RP}}\right) \leq
    \e_1 \right] \to 1$} as $n \to \infty$.
\end{theorem}

%
%
%
%
%
\noindent {\bf{Remarks.}}  Unlike the comparison against CQ, the correlation
matrix $R$ plays a large role in determining relative performance
of our test against SD.  Correlation enters in two different ways. First, 
the Frobenius norm $\mnorm{R}_F$ is larger when the data variables are more correlated. Second, if $\Sigma$
has a large number of small eigenvalues, then $\tr(D_{\sigma}^{-1})$
is very large when the variables are uncorrelated, i.e. when $\Sigma$ is
diagonal. Letting $U\Lambda U^{\top}$
be a spectral decomposition of $\Sigma$, with $u_i$ being the $i$th column of $U^{\top}$, note that $(D_{\sigma})_{ii} = u_{i}^{\top}\Lambda u_i$. When the data variables are correlated, the vector $u_i$ will have many nonzero components, which will give $(D_{\sigma})_{ii}$ a contribution from some of the larger eigenvalues of $\Sigma$, and prevent $(D_{\sigma})_{ii}$ from being too small. For example, if $u_i$ is uniformly distributed on the unit sphere, 
as in Example~\ref{random} below, 
then on average  $\E[(D_{\sigma})_{ii}] = \tr(\Sigma)/p$. Therefore, correlation has the effect of mitigating the
growth of $\tr(D_{\sigma}^{-1})$. Since the SD test
statistic \cite{SD2008} can be thought of as a version of the
Hotelling $T^2$ with a diagonal estimator of $\Sigma$, the SD test
statistic makes no essential use of correlation structure. By
contrast, our $T_k^2$ statistic \emph{does} take correlation into account,
 and so it is understandable that correlated data enhance the
performance of our test relative to SD.
%
%
%
%
\begin{example}\label{highCor}
Suppose the correlation matrix $R \in \R^{p\times p}$ has a block-diagonal structure, with
$m$ identical blocks $B\in \R^{d\times d}$ along the diagonal:
\begin{equation}\label{corrmatrix}
R=  \left( \begin{array}{ccc}
B & & \\
 & \ddots &  \\
 &      & B \end{array} \right).
\end{equation}
Note that $p=m \cdot d$. Fix a number $\rho \in (0,1)$, and let $B$ have
diagonal entries equal to 1, and off-diagonal entries equal to $\rho$,
i.e. $B = (1-\rho)I_{d\times d} + \rho {\bf{1}}{\bf{1}}^{\top}$, where ${\bf{1}}\in \R^d$
is the all-ones vector. Consequently, $R$ is positive-definite, and we may consider $\Sigma=R$
for simplicity. Since $\mnorm{B}_F^2 = d+ 2\rho^2 \binom{d}{2}$, and $\mnorm{R}_F^2 = m\mnorm{B}_F^2$, it follows that 
$$\mnorm{R}_F^2 = [1+\rho^2\,(d-1)] \, p.$$
Also, in this example we have $\tr(\Sigma)=\tr(D_{\sigma}^{-1})=p$ and $p/d =m$,  which implies
\begin{equation}
\begin{split}
\left(\frac{\tr(\Sigma)}{p}\right)^2 
\left(\frac{\tr(D_{\sigma}^{-1})} { \mnorm{R}_F} \right)^2 &= \frac{p}{1+\rho^2\,(d-1)}\leq \frac{m}{\rho^2}.
\end{split}
\end{equation} Under these conditions, we
conclude that the sufficient condition \eqref{suffSD} in Theorem
\ref{ThmSD} is satisfied when $n$ grows at a faster rate than the
number of blocks $m$.  Note too that the spectrum of $\Sigma$ consists
of $m$ copies of $\lambda_{\max}(\Sigma)=(1-\rho)+\rho\, d$ and
$(p-m)$ copies of $\lambda_{\min}(\Sigma)=1-\rho$, which means when
$\rho$ is not too small, the number of blocks is the same as the
number of dominant eigenvalues---revealing a parallel with Example
\ref{spiked}.  From these observations, it is straightforward to check
that this example satisfies assumptions ({\bf{A5}}) and ({\bf{A6}}) of Theorem~\ref{ThmSD}. The
simulations in Section \ref{sec:ROC} give an example where $R$ has the
form in line \eqref{corrmatrix} and the variables corresponding to
each block are highly correlated. 
\end{example}

\begin{example}\label{random} To consider the performance of our test in a case
where $\Sigma$ is not constructed deterministically,
Section~\ref{sec:ROC} illustrates simulations involving \emph{randomly
  selected} matrices $\Sigma$ for which $T_k^2$ is more powerful than
the tests of BS, CQ, and SD.  Random correlation
structure can be generated by sampling the matrix of eigenvectors of
$\Sigma$ from the uniform (Haar) distribution on the orthogonal group,
and then imposing various decay constraints on the eigenvalues of
$\Sigma$.  Additional details are provided in Section \ref{sec:ROC}.
\end{example}

\begin{example}\label{lowCor} It is possible to show that the sufficient
condition~\eqref{suffSD} \emph{requires} non-trivial correlation in
the high-dimensional setting. To see this, consider an example where the
data are completely free of correlation, i.e., where $R=I_{p\times
  p}$. Then, $\mnorm{R}_F=\sqrt{p}$, and Jensen's inequality implies that $\tr(D_{\sigma}^{-1})\geq
p^2/\tr(D_{\sigma})=p^2/\tr(\Sigma)$,  giving
${\small{\left(\frac{\tr(\Sigma)}{p}\right)^2 
\left(\frac{\tr(D_{\sigma}^{-1})} { \mnorm{R}_F} \right)^2}}\geq p$.
Altogether, this
shows if the data exhibits very low correlation, then \eqref{suffSD}
cannot hold when $p$ grows faster than $n$ (in the presence of a
uniformly oriented shift $\delta$). This is confirmed by the
simulations of Section~\ref{sec:ROC}. Similarly, it is shown in the
paper \cite{SD2008} that the SD test statistic is asymptotically
superior to the CQ test statistic\footnote{ Although the work in SD
  (2008) \cite{SD2008} was published prior to that of CQ (2010)
  \cite{CQ2010}, the asymptotic power function of CQ for problem \eqref{origproblem}
  is the same as
  that of the method proposed in BS (1996) \cite{BS96}, and SD offer a
  comparison against the method of BS.}  when $\Sigma$ is diagonal and
$\delta$ is a deterministic vector with all coordinates equal to the
same value.
\end{example}
%
%
%
The proof of Theorem~\ref{ThmSD} makes use of
concentration bounds for Gaussian quadratic forms, which are stated below in
Lemma~\ref{LemPoincare} (see Appendix~\ref{AppLemPoincare} for proof).
These bounds are similar to results in the
papers of Bechar, and Laurent and Massart~\cite{Bechar,LaurentMassart} (c.f. Lemma ~\ref{quadform} in Appendix 
\ref{ineqs}), but have error terms involving
the spectral norm as opposed to the Frobenius norm, and hence Lemma~\ref{LemPoincare}
may be of independent interest.

\begin{lemma}
\label{LemPoincare}
Let $A\in \R^{p\times p}$ be a positive semidefinite matrix with
$\mnorm{A}_2>0$, and let $Z\sim N(0,I_{p\times p})$. Then, for any
$t>0$,
\begin{equation}\label{lower}
\P\left[Z^{\top} A Z \geq \tr(A)\textstyle{\left(1+t\sqrt{\frac{\mnorm{A}_2}{\tr(A)}}\right)^2}\right] \leq \exp\left(-t^2/2\right),
\end{equation}
and for any $t \in \big (0, \sqrt{\frac{\tr(A)}{\mnorm{A}_2}-1}
\big)$, we have
\begin{equation}
\label{upper}
\P \left[Z^{\top} A Z \leq \tr(A) \textstyle{\left(\sqrt{ 1 - \frac{\mnorm{A}_2}
    {\tr(A)}} - t \sqrt{\frac{\mnorm{A}_2}{\tr(A)}} \right)^2 }\right]
\leq \exp(-t^2/2).
\end{equation}
\end{lemma}
%
%
\noindent Equipped with this lemma, we can now prove Theorem~\ref{ThmSD}.\\

\noindent \emph{Proof of Theorem \ref{ThmSD}.} We proceed along the lines of the proof of Theorem \ref{ThmCQ}.  Let us define the event
of interest, $\mathscr{E}_n \mydef \big \{
\text{ARE}\left(\beta_{\text{SD}}; \beta_{\text{RP}}\right) \leq \e_1
\big \}$, where we recall
$\text{ARE}\left(\beta_{\text{SD}}; \beta_{\text{RP}}\right)=
\Big({\frac{n
    \delta^{\top}D_{\sigma}^{-1}\delta}{\mnorm{R}_F}}\Big/\sqrt{n}
\Delta_k^2\Big)^2$ with $k=\lfloor n/2\rfloor$ and $y=1/2$.  The event $\mathscr{E}_n$ holds
if and only if
\begin{equation}\label{threeSecond}
\frac{n}{\mnorm{R}_F^2} \frac{1}{\e_1}\leq
\left(\frac{\Delta_k^2}{\|\delta\|_2^2}\right)^2\left(\frac{\|\delta\|_2^2}{\delta^{\top}D_{\sigma}^{-1}\delta}\right)^2.
\end{equation}
We consider the two factors on the right hand side of \eqref{threeSecond} separately. By Proposition~\ref{PropKeyLimit}, for any constant $c \in (0,
(1-\sqrt{a})^2)$, the first factor $\frac{\Delta_k^2}{\|\delta\|_2^2}$ satisfies
\begin{align}
\label{firstfactor}
\P\left( \frac{c\, k}{\tr(\Sigma)}
  \leq\frac{\Delta_k^2}{\|\delta\|_2^2} \right) & \to 1 \quad \mbox{
  as $n \rightarrow \infty$.}
\end{align}
Turning to the second factor
$\frac{\|\delta\|_2^2}{\delta^{\top}D_{\sigma}^{-1}\delta}$ in
line~\eqref{threeSecond}, we note that $\delta/\|\delta\|_2$
is uniformly distributed on the unit sphere of $\R^p$,
and so $\delta/\|\delta\|_2 \stackrel{d}{=}Z/\|Z\|_2$, where $Z\sim
N(0,I_{p\times p})$. Next, using Lemma~\ref{LemPoincare}, we see that
assumption ({\bf{A6}}) implies
$$\frac{Z^{\top} D_{\sigma}^{-1} Z}{\tr(D_{\sigma}^{-1})}\to 1 \text{ \ in probability.} $$
Since $\|Z\|_2^2/p\to 1$ almost surely, we obtain the limit
\begin{equation}\label{limit}
\frac{\delta^{\top}D_{\sigma}^{-1}\delta}{\|\delta\|_2^2}
\frac{p}{\tr(D_{\sigma}^{-1})}\stackrel{d}{=} \frac{Z^{\top}
  D_{\sigma}^{-1} Z}{\tr(D_{\sigma}^{-1})} \frac{p}{\|Z\|_2^2} \to 1
\text{ \ in probability.}
\end{equation}
Consequently, for any $\tilde{c}\in(0,1)$, the random variable 
$\frac{\|\delta\|_2^2}{\delta^{\top}D_{\sigma}^{-1}\delta}$ is greater than
$\frac{\tilde{c} \, p}{\tr(D_{\sigma}^{-1})}$ with probability tending to 1 as $n\to\infty$.
Applying this observation to line \eqref{threeSecond}, and using the limit \eqref{firstfactor}, 
we conclude that $\P\left(\mathscr{E}_n\right) \to
1$ as long as the inequality
\begin{equation}
\label{condition}
\frac{n}{\mnorm{R}_F^2} \frac{1}{\e_1} \leq \left(\frac{c\,
 k }{\tr(\Sigma)}\right)^2
\left(\frac{\tilde{c}\,p}{\tr(D_{\sigma}^{-1})}\right)^2
\end{equation}
holds for all large $n$. Replacing $k$ with $\frac{n}{2}\cdot \small{[1-o(1)]}$, 
the last condition is equivalent to 
\begin{equation}
 n\geq \left(\frac{\tr(\Sigma)}{p}\right)^2 
\left(\frac{\tr(D_{\sigma}^{-1})} { \mnorm{R}_F} \right)^2 \cdot \frac{4}{\e_1\,c^2 \, \tilde{c}^{\, 2} \, [1-o(1)]^2}.
\end{equation}
Thus, for a given choice of $c_1(\e_1)$ in the statement 
of the theorem, it is possible to choose  $c
< (1-\sqrt{a})^2$ and $\tilde{c}<1$ so that the claimed sufficient condition~\eqref{suffSD} implies
the inequality~\eqref{condition} for all large $n$, which completes
the proof.  \qed

%
%
%
%
\section{Performance comparisons on real and synthetic data}\label{sec:data}
In this section, we compare our procedure to a broad collection of competing methods
 on synthetic data, illustrating the effects of the different factors involved in
Theorems~\ref{ThmCQ} and \ref{ThmSD}. Sections \ref{sec:ROC}  and \ref{sec:calibration} 
consider ROC curves and calibration curves respectively. An example involving high-dimensional 
gene expression data is studied in Section \ref{sec:gene}.

\subsection{ROC curves on synthetic data}
\label{sec:ROC}
\begin{figure}[h!]
\begin{center}
\includegraphics[angle=0,width=0.5\linewidth]{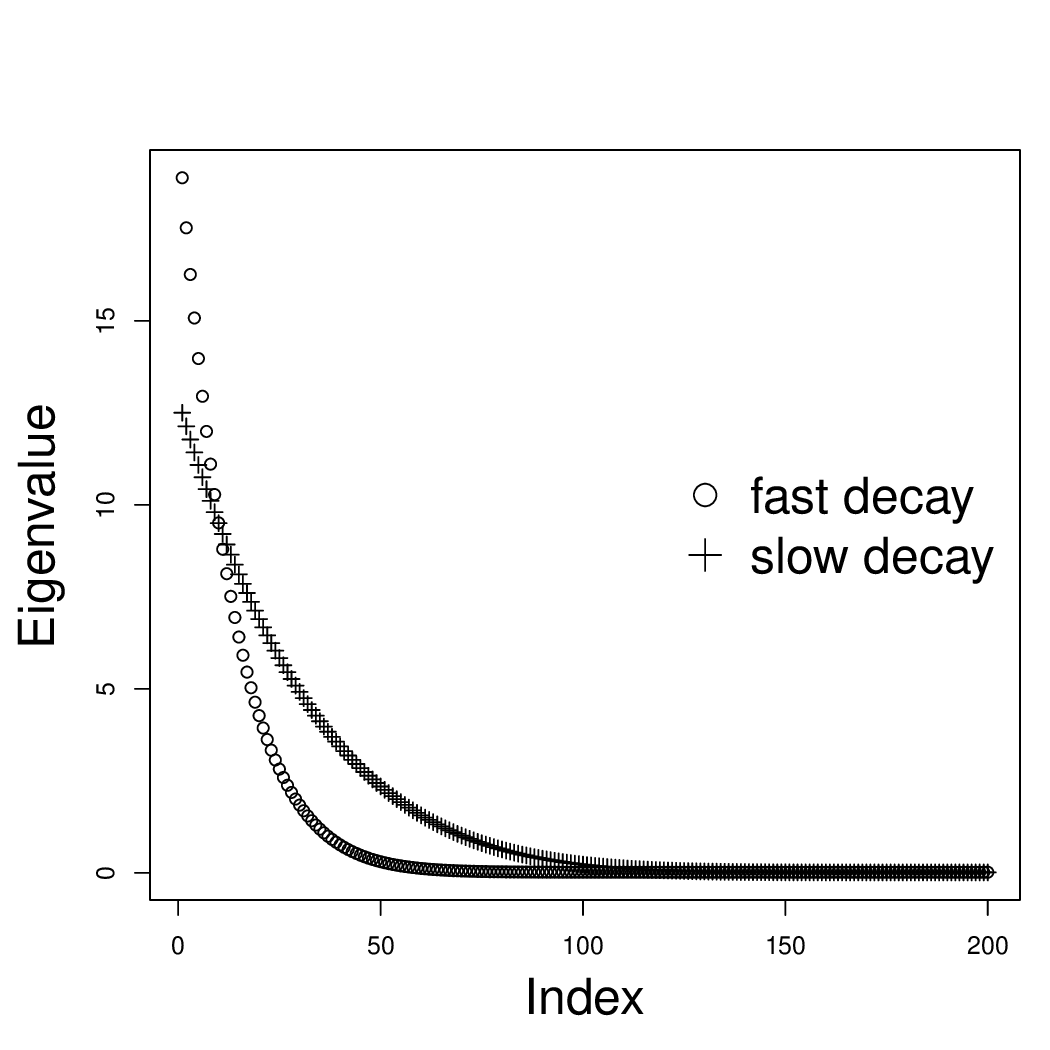}
\caption{Plots of two sets of eigenvalues $\lambda_1,\dots,\lambda_p$,
  with slow and fast decay, both satisfying
  $\mnorm{\Sigma}_F=\sqrt{\lambda_1^2+\cdots+\lambda_p^2}=50$. To interpret the number
  of non-negligible eigenvalues, there are
29 eigenvalues greater than $\frac{1}{10}\lambda_{\max}(\Sigma)$ in the case of fast decay,
and there are 65 eigenvalues greater than $\frac{1}{10}\lambda_{\max}(\Sigma)$ in the case of 
slow decay.}
\label{fig:spectra}
\end{center}
\end{figure}

\begin{figure}[!]
\centering
\subfloat[font=Large][diagonal $\Sigma$, slow decay]{\includegraphics[angle=-90,
  width=.35\linewidth]{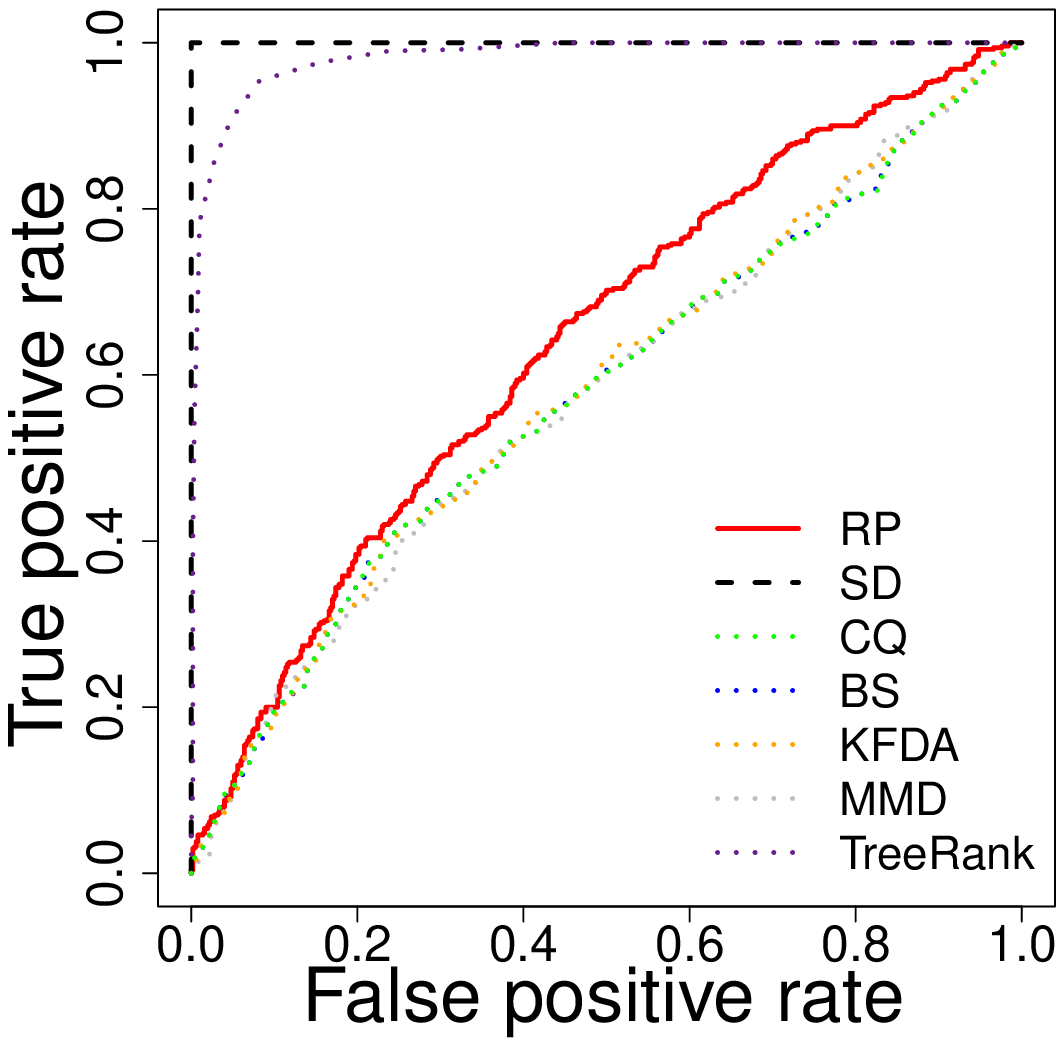}} \ \ \ \ \ 
\subfloat[diagonal $\Sigma$, fast decay]{\includegraphics[angle=-90,
  width=.35\linewidth]{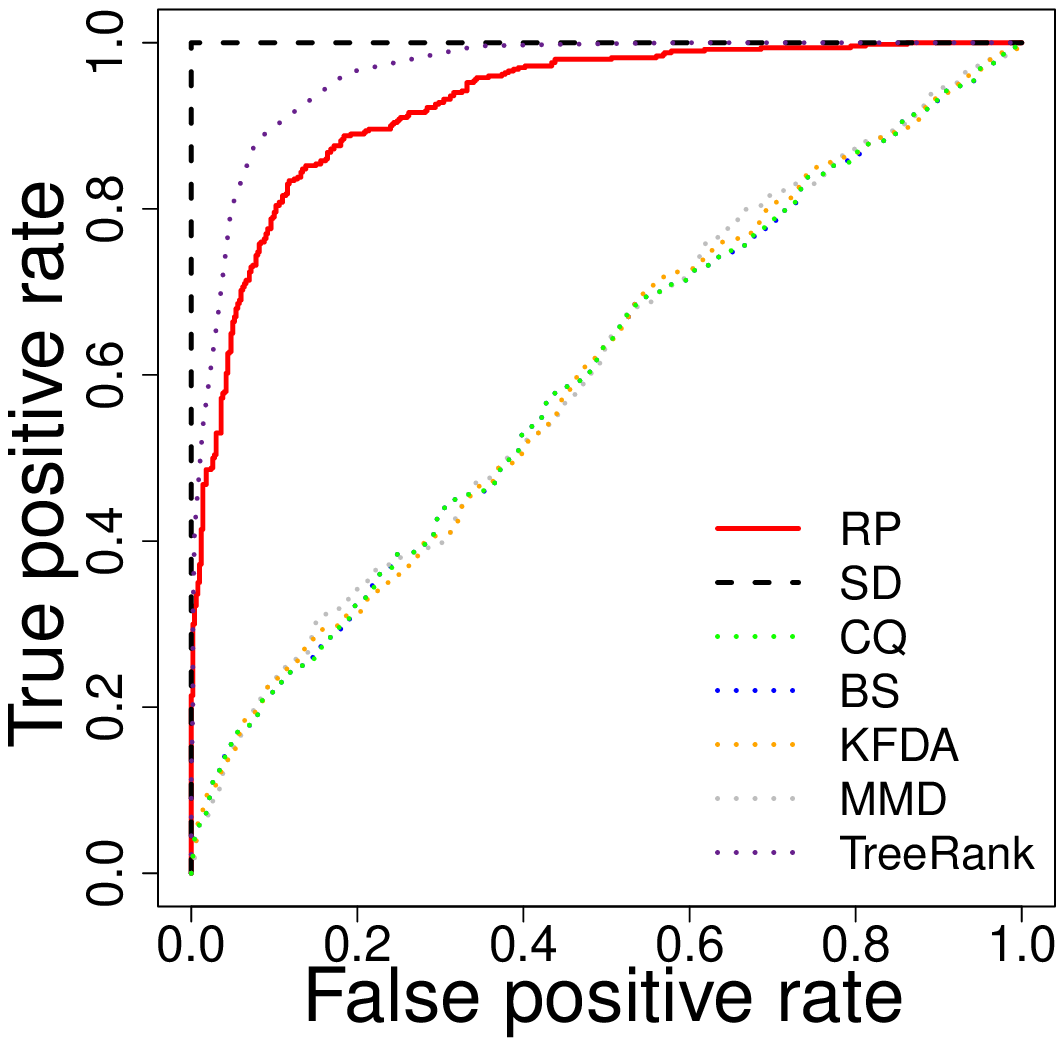}}\\
\subfloat[random $\Sigma$, slow decay]{\includegraphics[angle=-90, width=.35\linewidth]{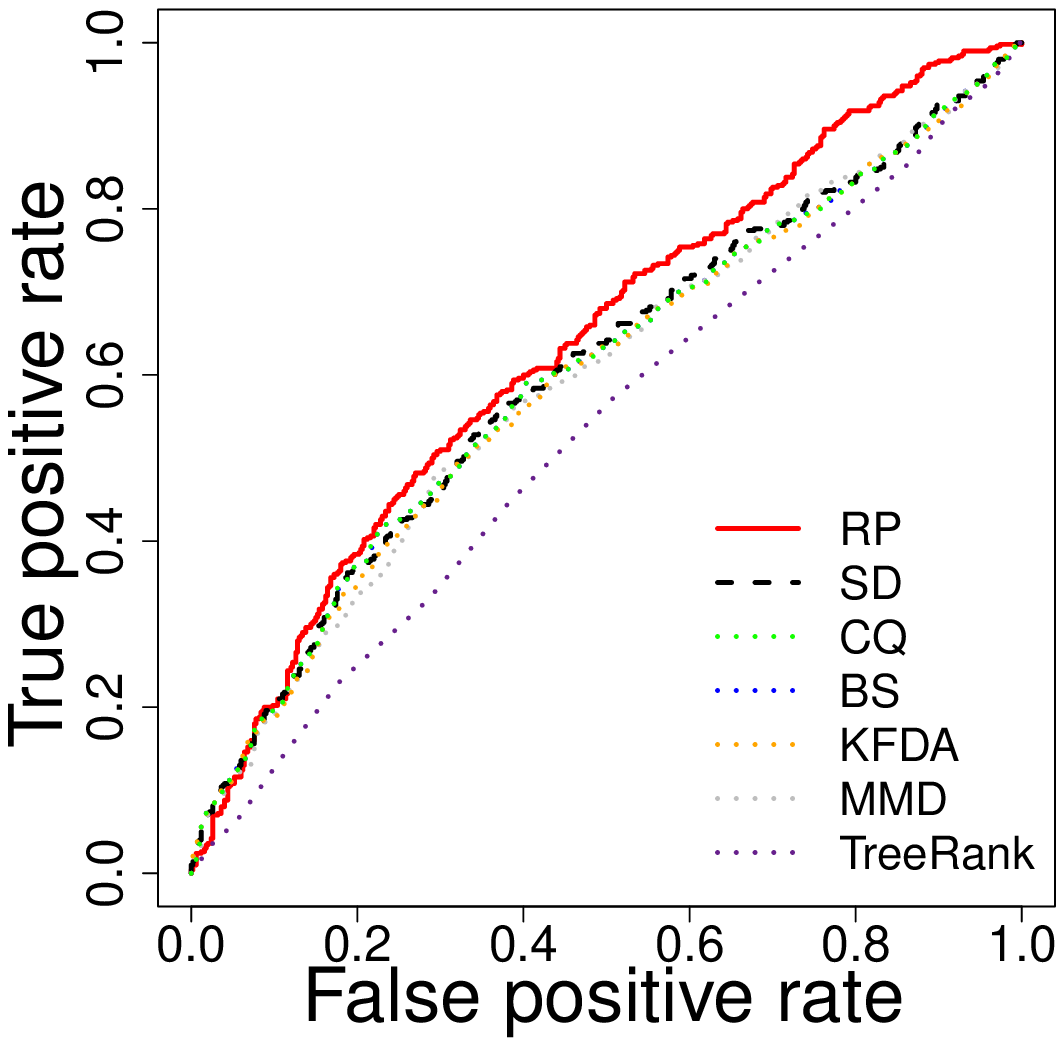}} \ \ \ \ 
\subfloat[random $\Sigma$, fast decay]{\includegraphics[angle=-90, width=.35\linewidth]{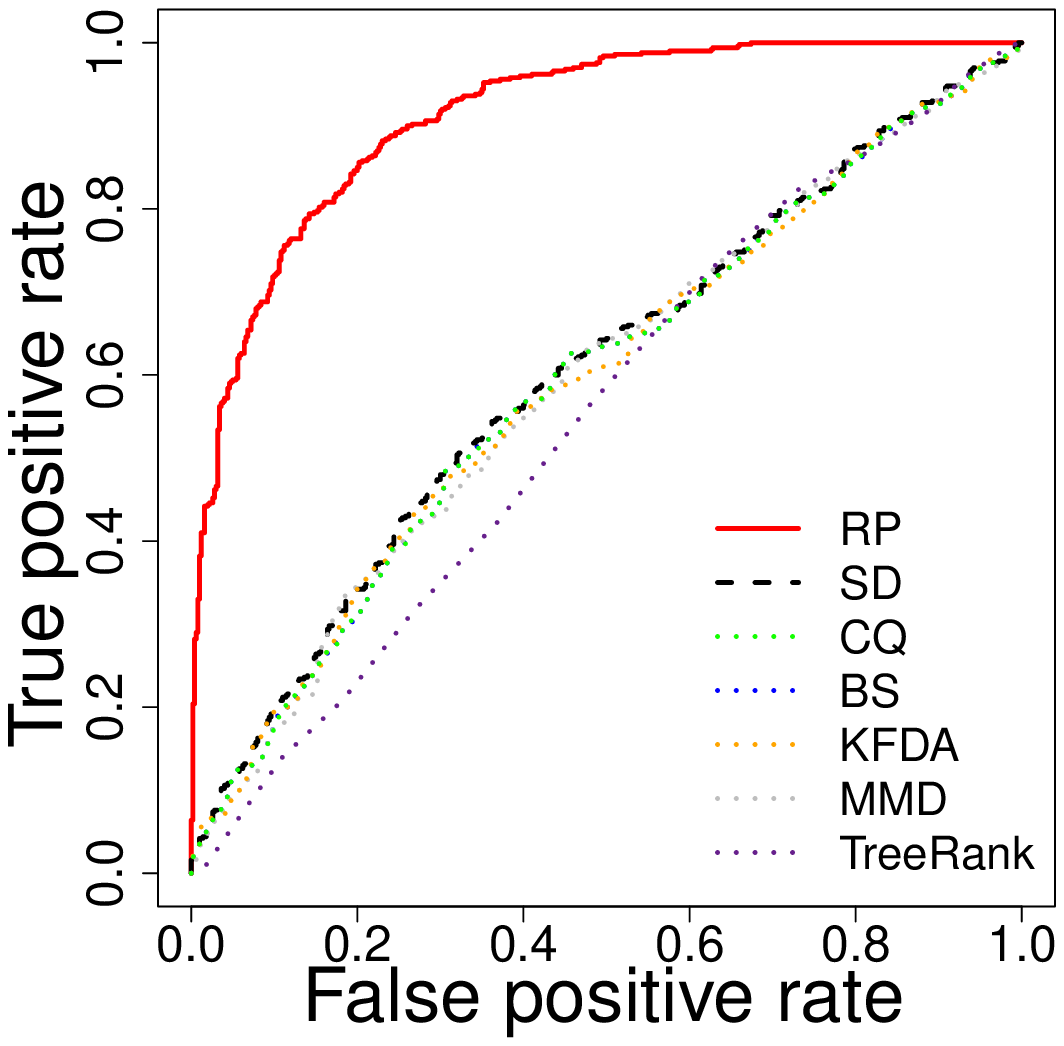}}\\
  \subfloat[block-diagonal correlation]{\includegraphics[angle=-90, width=0.35\linewidth]{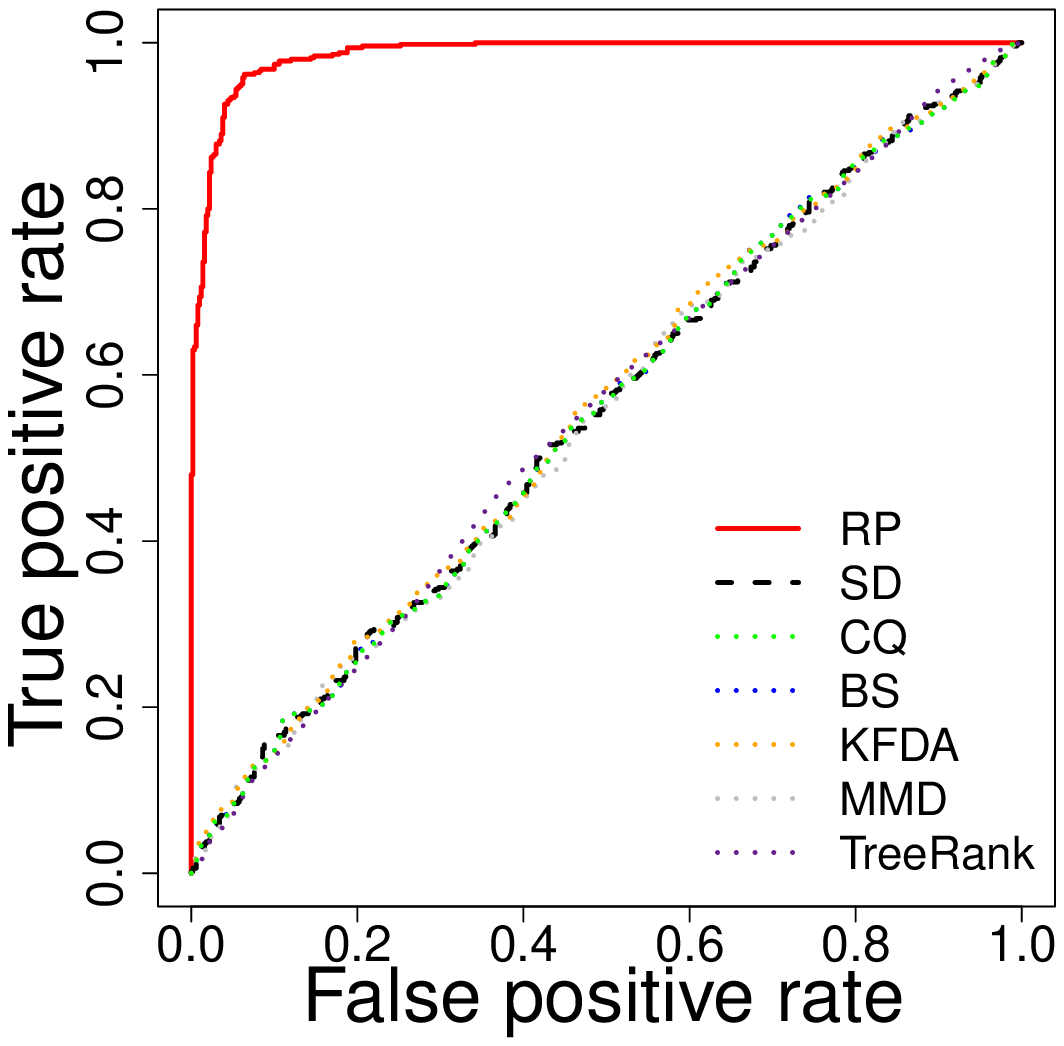}}
\caption{ ROC curves of several test statistics
  for five different settings of correlation structure and spectral decay of $\Sigma$: (a) Diagonal covariance / slow decay, (b) Diagonal covariance /
  fast decay, (c) Random covariance / slow decay, and (d) Random covariance /
  fast decay. (e) Block-diagonal correlation.}
\label{fig:exps}
\end{figure}

%
Using multivariate normal data, we generated ROC curves (see Figure \ref{fig:exps}) in 
five distinct parameter settings.
 For each ROC curve, we sampled $n_1=n_2=50$ 
data points from each of the 
distributions $N(\mu_1,\Sigma)$ and $N(\mu_2,\Sigma)$ in $p=200$ dimensions, 
and repeated the process $500$ times with $\delta=\mu_1-\mu_2=0$ under $\hh_0$, and $500$ times
with $\|\delta\|_2=1$ under $\hh_1$. For each simulation under $\hh_1$, the shift $\delta$ was sampled
as $Z/\|Z\|_2$ for $Z\sim N(0,I_{p\times p})$, so as to be drawn uniformly from the unit sphere, and
satisfy assumption ({\bf{A3}}) in Theorems~\ref{ThmCQ}
and \ref{ThmSD}. Letting $U\Lambda U^{\top}$ denote a spectral decomposition of $\Sigma$, 
we specified the first four parameter settings by choosing $\Lambda$ to have a 
spectrum with slow or fast decay, and choosing $U$ to be $I_{p\times p}$ or
a randomly drawn matrix from the uniform (Haar) distribution on the orthogonal group
\cite{Stewart1980Efficient}. Note that $U=I_{p\times p}$ gives a diagonal covariance matrix
$\Sigma$, whereas a randomly chosen $U$ induces correlation among the variables. 
To consider two rates of spectral decay, we selected $p$ equally spaced eigenvalues
$\lambda_1,\dots,\lambda_p$ between $10^{-2}$ and $1$, and raised them
to the power $15$ for fast decay, and the power $6$ for slow decay. We
then added $10^{-3}$ to each eigenvalue to control the condition
number of $\Sigma$, and rescaled them so that $\mnorm{\Sigma}_F=\sqrt{\lambda_1^2+\cdots+\lambda_p^2}=50$ in each of the first four settings (fixing a common amount of variance).
Plots of the resulting spectra are shown in Figure \ref{fig:spectra}.
The fifth setting was specified by choosing
the correlation matrix $R$ to have a block-diagonal structure, corresponding
to 40 groups of highly correlated variables. Specifically, the matrix
$R$ was constructed to have 40 identical blocks $B\in \R^{5\times 5}$ along its diagonal,
with the diagonal entries of $B$ equal to 1, and the off-diagonal entries of $B$ equal to $\rho:=\frac{1}{1.01}$ (c.f. Example \ref{highCor}). The matrix $\Sigma$ was then formed by setting 
$D_{\sigma}= \frac{1}{\rho}\, I_{p\times p}$, and $\Sigma = D_{\sigma}^{1/2} R D_{\sigma}^{1/2}$.

In addition to our random projection (RP)-based test, we implemented
the methods of BS \citep{BS96}, SD \citep{SD2008}, and CQ
\citep{CQ2010}, which are all designed specifically for problem
\eqref{origproblem} in the high-dimensional setting. For the sake of
completeness, we also show comparisons against two recent
non-parametric procedures that are based on kernel methods: maximum
mean discrepancy (MMD) \citep{Gretton2007A}, and kernel Fisher
discriminant analysis (KFDA) \citep{Harchaoui2008Testing}, as well as
a test based on area-under-curve maximization, denoted TreeRank
\citep{AUCopt}. Overall, the ROC curves in Figure~\ref{fig:exps} show that in each of the five settings, either our test, or the test of SD, perform the best within this collection of procedures.

On a qualitative level, Figure~\ref{fig:exps} reveals some striking
differences between our procedure and the competing tests.  Comparing
independent variables versus correlated variables, i.e.\:panels (a) and
(b), with panels (c) and (d), we see that the tests of SD and TreeRank
lose power in the presence of correlated data. Meanwhile, the ROC
curve of our test is essentially unchanged when passing from
independent variables to correlated variables. Similarly, our test
also exhibits a large advantage when the correlation structure is
prescribed in a block-diagonal manner in panel (e). The agreement of
this effect with Theorem~\ref{ThmSD} is explained in the remarks and
examples after that theorem. Comparing slow spectral decay versus fast
spectral decay, i.e.\:panels (a) and (c), with panels (b) and (d), we
see that the competing tests are essentially insensitive to the change
in spectrum, whereas our test is able to take advantage of
low-dimensional covariance structure. The remarks and examples of
Theorem~\ref{ThmCQ} give a theoretical justification for this
observation.

It is also possible to offer a more quantitative assessment of the ROC curves
in light of Theorems ~\ref{ThmCQ} and~\ref{ThmSD}. Table~\ref{tab:are} summarizes approximate values of $\tr(\Sigma)^2/\mnorm{\Sigma}_F^2$ and {\footnotesize{$\left(\frac{\tr(\Sigma)}{p}\right)^2\left(\frac{\tr(D_{\sigma}^{-1})}{\mnorm{R}_F}\right)^2$}} from Theorems~\ref{ThmCQ} and~\ref{ThmSD}
in the five settings described above.\footnote{For the case of randomly selected $\Sigma$, the quantities are obtained as the average from 500 draws.} The table shows that our theory is consistent with Figure \ref{fig:exps}\, in the sense that the only settings for which our test yields an inferior ROC curve are those for which the quantity {\footnotesize{$\left(\frac{\tr(\Sigma)}{p}\right)^2\left(\frac{\tr(D_{\sigma}^{-1})}{\mnorm{R}_F}\right)^2$}} is drastically larger than $n=50+50-2=98$. (In all of the settings where $\tr(\Sigma)^2/\mnorm{\Sigma}_F^2$ and {\footnotesize{$\left(\frac{\tr(\Sigma)}{p}\right)^2\left(\frac{\tr(D_{\sigma}^{-1})}{\mnorm{R}_F}\right)^2$}} are less than $n$, our test yields the best ROC curve against the competitors.) However, if the entries in the table are multiplied by a choice of the constant $c_1(\e_1)>\frac{4}{\e_1(1-\sqrt{a})^4}$ from Theorems~\ref{ThmCQ} and~\ref{ThmSD}, we see that our asymptotic conditions \eqref{suffCQ} and \eqref{suffSD} are somewhat conservative at the finite sample level. Considering $n=98$, the table shows that $c_1(\e_1)$ would need to be roughly equal to 1.5 so that the inequalities \eqref{suffCQ} and \eqref{suffSD} hold in all the settings for which our method has a better ROC curve than the relevant competitor. In the basic case that $\e_1=1$ and $a=0$, we have $\frac{4}{\e_1(1-\sqrt{a})^4}=4$, which means that the constant $\frac{4}{\e_1(1-\sqrt{a})^4}$ needs to be improved by roughly a factor of $4/1.5\simeq 3$ or better. We expect that such improvement is possible with a more refined analysis of the proof of Proposition~\ref{PropKeyLimit}.
\begin{table}[h!]
\footnotesize
\begin{center}
\caption{Approximate values of the quantities
  {\footnotesize{$\tr(\Sigma)^2/\mnorm{\Sigma}_F^2$}} and
  {\footnotesize{$\left(\frac{\tr(\Sigma)}{p}\right)^2\left(\frac{\tr(D_{\sigma}^{-1})}{\mnorm{R}_F}\right)^2$}}
  in the five parameter settings of the synthetic data
  experiments. Theorems~\ref{ThmCQ} and \ref{ThmSD} assert that these
  quantities determine the relative performance of our test against CQ
  and SD respectively.}
\label{tab:are}
  \vspace{10pt}
\renewcommand{\arraystretch}{1.2}

\begin{tabular}{| l | c | c | c | c | c | }
\hline
                         & diagonal $\Sigma$,  &  diagonal $\Sigma$,  &  random $\Sigma$,  & random $\Sigma$, & block-diagonal \\
                         & slow decay & fast decay & slow decay & fast decay & correlation \\ \hline
 (Thm.~\ref{ThmCQ} vs. CQ) \, \footnotesize{$\tr(\Sigma)^2/\mnorm{\Sigma}_F^2$}& 54 & 25 & 54 & 25 & 41\\
 \hline
(Thm.~\ref{ThmSD} vs. SD) \footnotesize{$\left(\frac{\tr(\Sigma)}{p}\right)^2\left(\frac{\tr(D_{\sigma}^{-1})}{\mnorm{R}_F}\right)^2$ }           & $4.6\times 10^5$ &$3.5\times 10^5$ & 58 & 30 &  41 \\
\hline
\end{tabular}
\end{center}
\end{table}
\subsection{Calibration curves on synthetic data}\label{sec:calibration}
\begin{figure}[h]
\centering \subfloat[font=Large][diagonal $\Sigma$, slow
  decay]{\includegraphics[angle=-90,
    width=.3\linewidth]{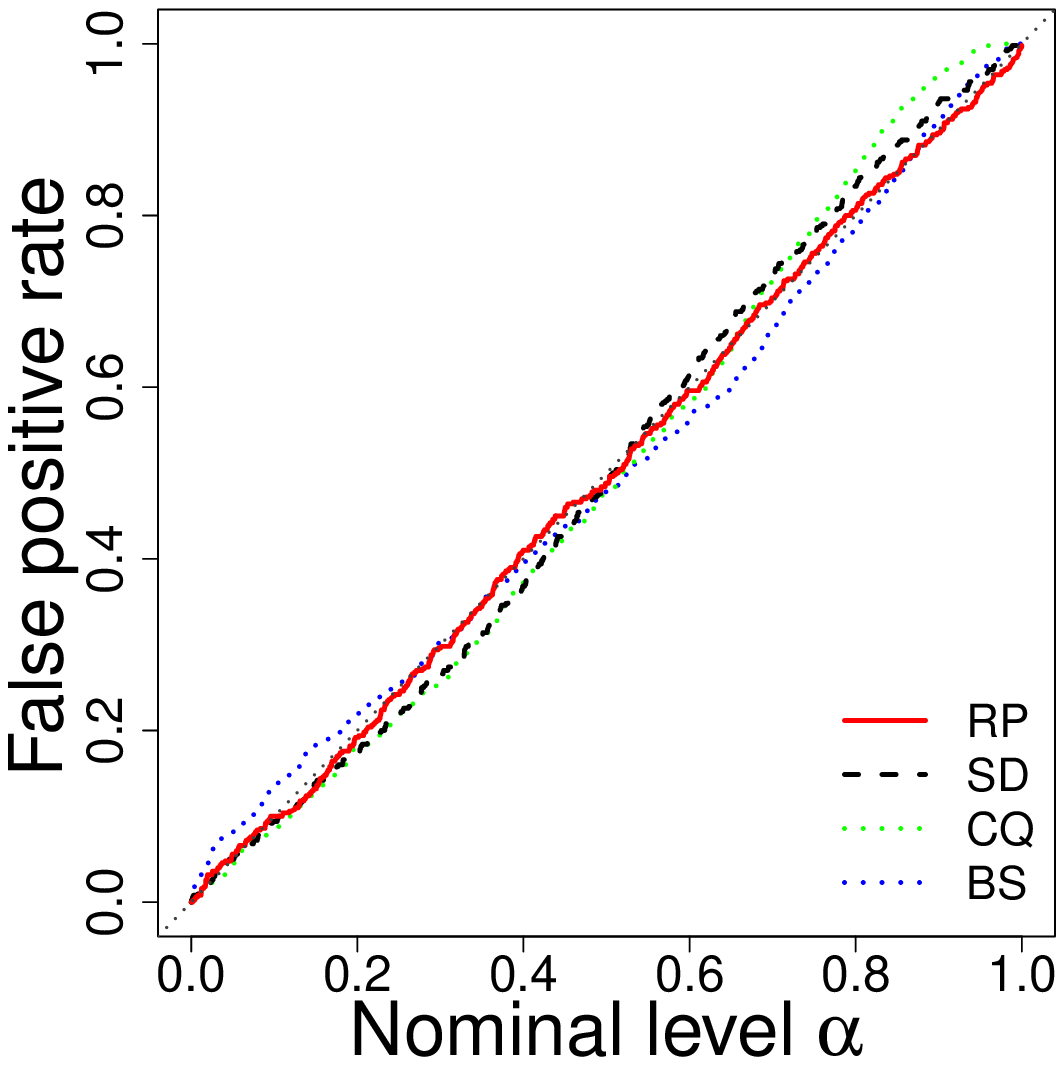}}
\ \ \ \ \subfloat[font=Large][random $\Sigma$, fast
  decay]{\includegraphics[angle=-90,
    width=.3\linewidth]{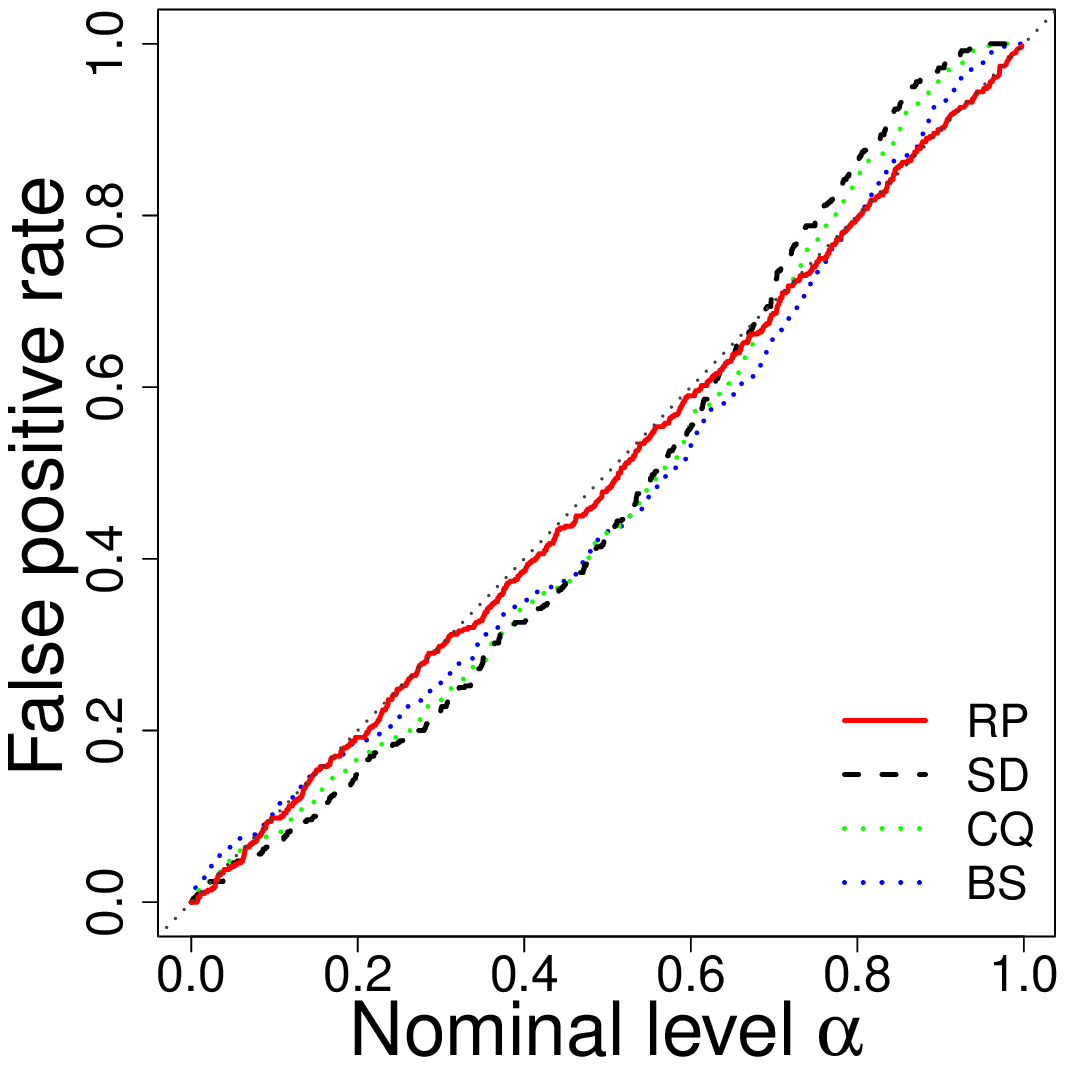}}\ \ \ \ \subfloat[font=Large][mixture
  model]{\includegraphics[angle=-90,width=.3\linewidth]{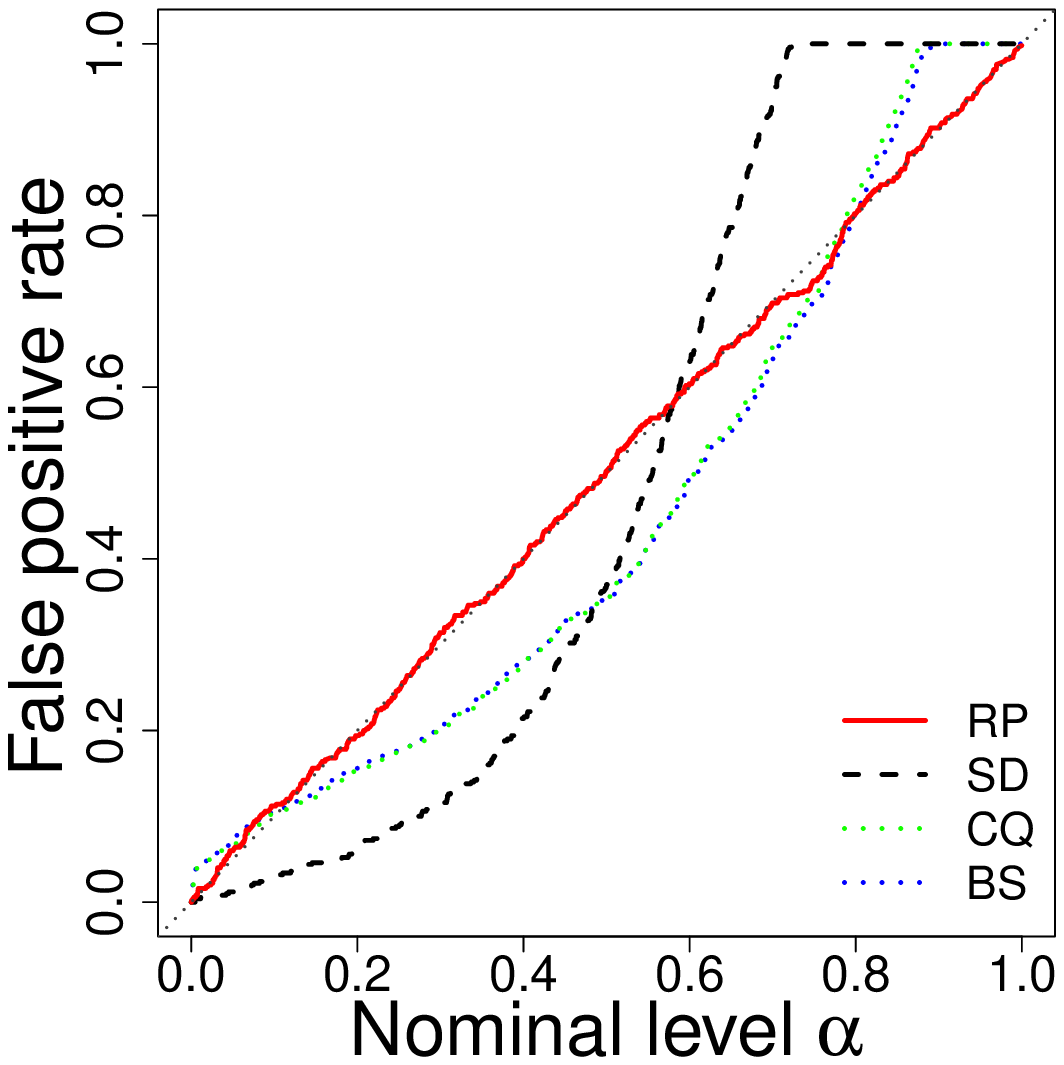}}
  \vspace{5pt.}
\caption{Calibration plots on synthetic data under three different data-generating distributions. The grey line is a reference for optimal calibration.}
\label{fig:synCalibration}
\end{figure}

Figure~\ref{fig:synCalibration}\:contains calibration plots resulting from the simulations
described in Section \ref{sec:ROC}---showing how well the observed
false positive rates (FPR) of the various tests compare against the nominal level $\alpha$.
(Note that these plots only reflect simulations under $\hh_0$.)
Ideally, when testing at level $\alpha$, the observed
FPR should be as close to $\alpha$ as possible,
and a thin diagonal
grey line is used here as a reference for perfect calibration.
Figures~\ref{fig:synCalibration} (a) and (b) correspond respectively to the settings from Section~\ref{sec:ROC} where $\Sigma$ is diagonal, with a slowly decaying spectrum, and where $\Sigma$ has
random eigenvectors and a rapidly decaying spectrum. In these cases
the tests of BS, CQ, and SD are reasonably well-calibrated, and our
test is nearly on top of the optimal diagonal line. To consider
robustness of calibration, we repeated the simulation from panel (a), but replaced the sampling distributions $N(\mu_i,\Sigma)$, $i=1,2$, with the mixtures \mbox{$0.2\, N(\mu_i+d_{1,i},\Sigma)+0.3\, N(\mu_i+d_{2,i},\Sigma)+ 0.5 \, N(\mu_i+d_{3,i},\Sigma)$}, where $0.2\,d_{1,i}+0.3\,d_{2,i}+0.5\,d_{3,i}=0$, and $d_{1,1},d_{2,1},d_{1,2},d_{2,2}$ were drawn independently and uniformly from a sphere of radius $\mnorm{\Sigma}_F^2$.
The resulting calibration plot in Figure~\ref{fig:synCalibration} (c) shows that our test deviates slightly from the diagonal in this case, but the calibration of the other three tests degrades to a much more noticeable extent. Experiments on other non-Gaussian distributions (e.g. with heavy tails) gave similar results, 
suggesting that the critical values of our procedure may be generally more robust (see also the discussions
of robustness in Sections \ref{sec:intuition} and \ref{sec:gene}).
%
%
%
%
\subsection{Comparison on high-dimensional gene expression data.}
\label{sec:gene}
The ability to detect gene sets having different expression between
two types of conditions, e.g., benign and malignant forms of a
disease, is of great value in many areas of biomedical research. 
In this section, we study our testing procedure in the context of
determining whether a set of $p$ genes is differentially expressed
between two relatively small groups of patients of sizes $n_1$ and
$n_2$.
To compare the performance of our $T_k^2$ statistic against competitors
CQ and SD in this type of application, we constructed a
collection of 1680 distinct two-sample problems in the following
manner, using data from three genomic studies of
ovarian~\citep{Tothill2008Novel}, myeloma~\citep{Moreaux2011high-risk}
and colorectal~\citep{Jorissen2009Metastasis-Associated} cancers. 
First, we randomly split the 3 datasets respectively into $6$, $4$,
and $6$ groups of approximately $50$ patients. Next, we considered all
possible pairwise comparisons between all sets of patients on each of 14
biologically meaningful gene sets from the canonical pathways of the database
MSigDB~\citep{Subramanian2005Gene}.  Each gene set contains between 
$75$ and $128$ genes (with an average of $98.5$). Since $n_1\simeq n_2\simeq 50$, our
collection of two-sample problems is genuinely high-dimensional.
Specifically, we have $14\times
(\binom{6}{2}+\binom{4}{2}+\binom{6}{2}) = 504$ problems under $\hh_0$, and $14\times(6\times 4+6\times4+6\times6)= 1176$ problems
under $\hh_1$, where we assume that every gene set is differentially
expressed between two sets of patients with two different cancers, and
that no gene set is differentially expressed between two sets of
patients with the same cancer. Although it is conceivable that this assumption
could be violated by the existence of various cancer subtypes, or differences between original
tissue samples, our initial step of randomly splitting the three cancer datasets into subsets guards against this possibility.

\begin{figure}[!]
\centering
\subfloat[FPR for genomic data]{\includegraphics[angle=270, width=.4\linewidth]{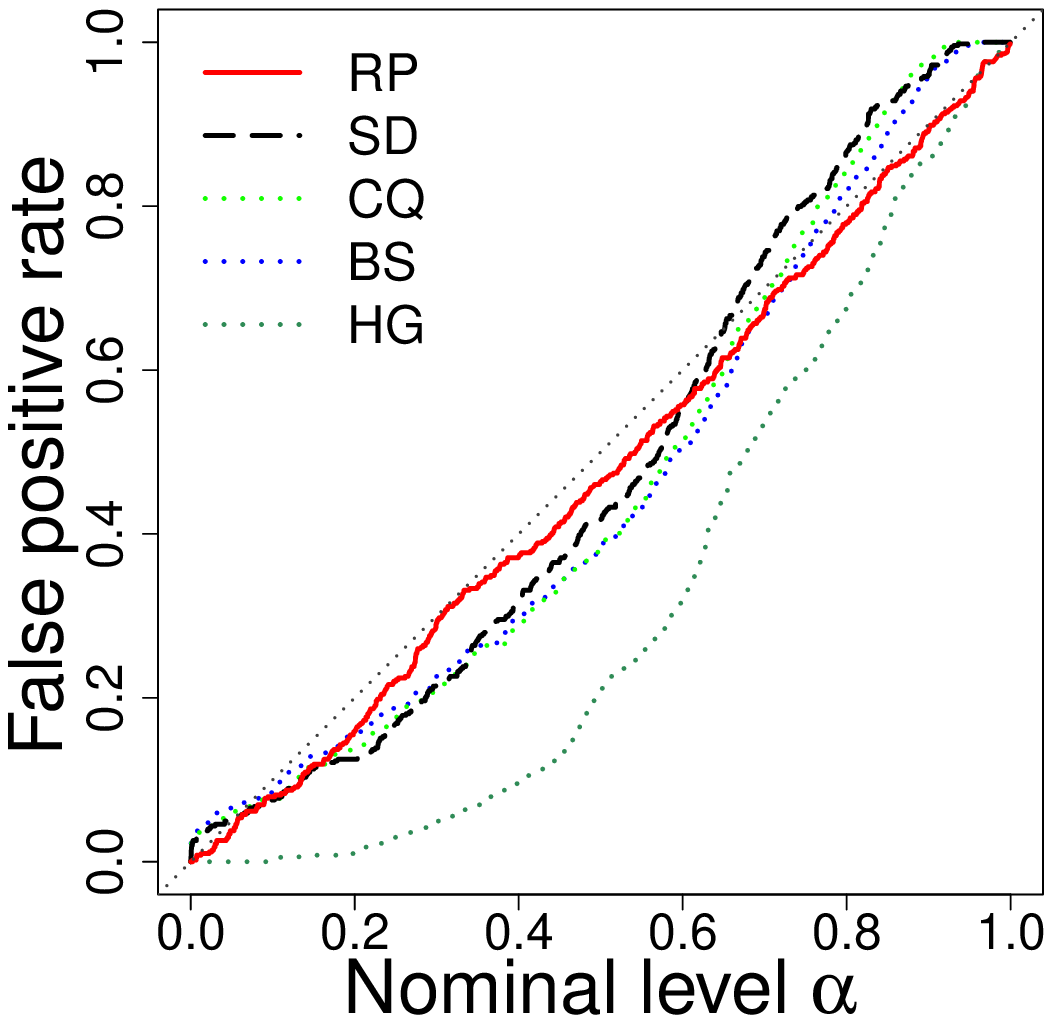}} \ \ \ 
\subfloat[FPR for genomic data (zoom)]{\includegraphics[angle=270,
  width=.4\linewidth]{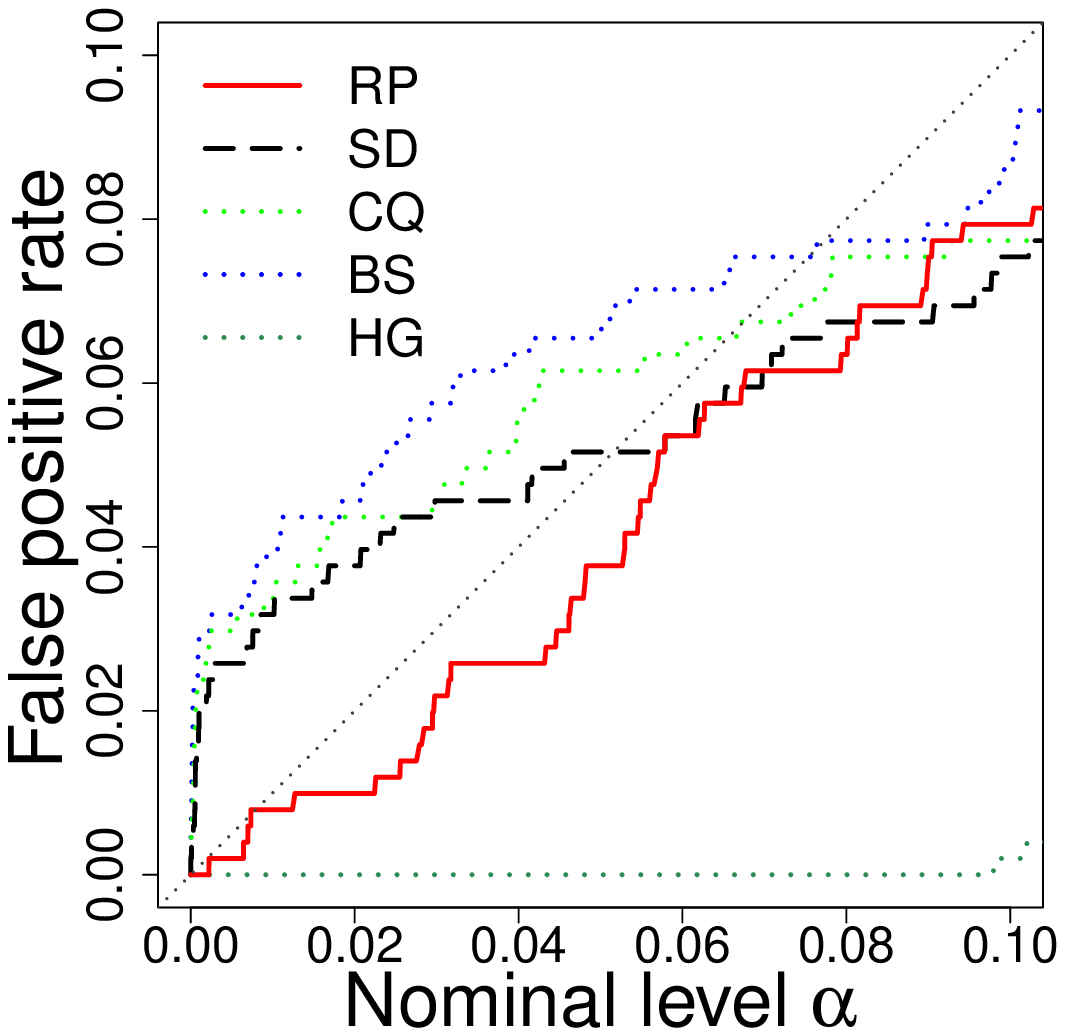}}
  \caption{(a) False positive rate against p-value
  threshold on the gene expression experiment of Section \ref{sec:gene} for
  RP, BS, CQ, SD and the hypergeometric (HG) enrichment test. (b) Zoom on the
  $\textrm{p-value}<0.1$ region. The optimal diagonal line is plotted in light grey.}
\label{fig:realdata}
\end{figure}

With consideration to ROC curves, the
cancer datasets are dissimilar enough that BS, CQ, SD, and our method all
produce perfect ROC curves from the collection of two-sample problems 
(no $\hh_1$ case has a larger p-value than any $\hh_0$ case). The 
hypergeometric test-based (HG) enrichment analysis \citep{Beissbarth2004GOstat} 
often used by experimentalists on this problem gives a suboptimal
area-under-curve of $0.989$.

Examining the quality of calibration reveals an important difference between our
test and the competitors in this example.
 It is apparent in Figure~\ref{fig:realdata}~(a) that the curve
for our procedure is closer to the
optimal diagonal line (plotted in light grey) for most values of $\alpha$ than the competing curves.  
 Furthermore, the lower-left corner of Figure
~\ref{fig:realdata}~(a) is of particular importance, as practitioners are
usually only interested in p-values lower than $10^{-1}$.
Figure~\ref{fig:realdata}~(b) is a zoomed plot of the lower-left corner,
which shows that the SD and CQ tests commit too many false positives
at low thresholds. Again, in this regime, our procedure is closer to
the diagonal and safely commits fewer than the allowed number of false positives. 
For example, when thresholding p-values at $0.01$, SD
has an actual FPR of $0.03$, and an even more excessive FPR of $0.02$ 
when thresholding at $0.001$. The tests of CQ and BS do even worse. 
The same thresholds on the p-values of our test lead to false positive 
rates of $0.008$ and $0$ respectively.

As discussed in Section \ref{sec:intuition}, there are two properties of our testing procedure that could account
for the advantage of our FPR on the both the synthetic and real data. 
First, our test inherits exact critical values for Gaussian data from the classical
Hotelling test, whereas the competing tests of SD, CQ, and BS use
thresholds based on asymptotic approximations. Second, even if the
$p$-dimensional data is poorly approximated by $N(\mu_1,\Sigma)$ and
$N(\mu_1,\Sigma)$, it is well known that randomly projected data tends
to be nearly Gaussian \cite{FreedmanDiaconis}. Consequently, the use
of a projection that induces Gaussianity, in conjunction with exact
critical values for Gaussian data may explain the advantage of our
test's FPR.
%
%
%
%
%
%
%
\section{Conclusion}
We have proposed a novel testing procedure for the two-sample test of
means in high dimensions. This procedure can be implemented in a
simple manner by first projecting a dataset with a single randomly
drawn matrix, and then applying the standard Hotelling $T^2$ test in
the projected space.  In addition to deriving an asymptotic power function for
this test, we have provided interpretable conditions on the covariance
and correlation matrices for achieving greater power
than competing tests in the sense of asymptotic relative efficiency.
Specifically, our theoretical comparisons show that our test is well-suited 
to interesting regimes where the data variables are correlated, or where 
most of the variance can be captured in a small number of variables. 
Furthermore, in the realistic case of $(n,p)=(100,200)$, these types of
conditions were shown to correspond to favorable performance of our test against
several competitors in ROC curve comparisons on synthetic
data. Finally, we showed on real gene expression data that our
procedure was more reliable than competitors in terms of its false
positive rate.  Extensions of this work may include more refined
applications of random projection to other high-dimensional testing
problems.
\paragraph{Acknowledgements.}  The authors 
thank Sandrine Dudoit, Anne Biton, and Peter Bickel for helpful
discussions. MEL gratefully acknowledges the support of the DOE CSGF
Fellowship, under grant number DE-FG02-97ER25308, and LJJ the support
of Stand Up to Cancer. MJW was partially supported by NSF
grant DMS-0907632.

\newpage
\appendix
\section{Matrix and Concentration Inequalities}
\label{ineqs}

This appendix is devoted to a number of  matrix and
concentration inequalities used at various points in our analysis.
We also prove Lemma \ref{LemPoincare}, which is stated in the main text in Section~\ref{sec: sdcomp}.

\begin{lemma}\label{traceprod}
If $A$ and $B$ are square real matrices of the same size with $A\succeq
0$ and $B = B^{\top}$, then
\begin{equation}
\lambda_{\min}(B) \tr(A) \leq \tr(AB) \leq \lambda_{\max}(B) \tr(A).
\end{equation}
\end{lemma}
\noindent {\emph{Proof.}} The upper bound is an immediate consequence
of Fan's inequality \cite[p.10]{BorweinLewis}, which states that any
two symmetric matrices $A,B\in \R^{p\times p}$ satisfy $\tr(AB)\leq
\sum_{i=1}^p \lambda_i(A)\lambda_i(B),$ where $\lambda_{i+1}(\cdot)\leq \lambda_i(\cdot)$.
Replacing $A$ with $-A$ yields the lower bound. \qed\\

%
%
%

\noindent See the papers of Bechar~\cite{Bechar}, or Laurent and
Massart~\cite{LaurentMassart} for proofs of the following
concentration bounds for Gaussian quadratic forms.

\begin{lemma}\label{quadform} Let $A \in \R^{p\times p}$ with $A\succeq 0$, 
and $Z\sim N(0,I_{p\times p})$. Then for any $t > 0$, we have
\begin{subequations}
\begin{align}
\P \left[ Z^{\top}AZ \geq \tr(A) + 2\mnorm{A}_F \sqrt{t} +
  2\mnorm{A}_2 t \right] & \; \leq \; \exp(-t), \quad \mbox{and} \\
\P \left[ Z^{\top}AZ \leq \tr(A)-2\mnorm{A}_F\sqrt{t} \right] & \;
\leq \; \exp(-t).
\end{align}
\end{subequations}
\end{lemma}
%
%
\noindent The following result on the extreme eigenvalues of Wishart
matrices is given in Davidson and Szarek \cite[Theorem
  II.13]{DavidsonSzarek}.

\begin{lemma}[]
\label{wishart} For $k \leq p$, let $P_k^{\top} \in \R^{k\times p}$ be a random 
matrix with i.i.d. $N(0,1)$ entries. Then, for all $t > 0$, we have
\begin{subequations}
\begin{align}
\textstyle
 \P \left[ \lambda_{\max} \big( \frac{1}{p}P_k^{\top}P_k\right) \geq
   \left(1+\sqrt{k/p} +t \big)^2 \right] & \leq \exp(-pt^2/2), \quad
 \mbox{and}\\
\textstyle
\P \left[ \lambda_{\min} \big( \frac{1}{p} P_k^{\top} P_k \big) \leq
  \big(1- \sqrt{k/p} - t \big)^2 \right] & \leq \exp(-pt^2/2).
\end{align}
\end{subequations}
\end{lemma}
~\\
\noindent{\bf{Proof of Lemma~\ref{LemPoincare}.}}\label{AppLemPoincare}
Note that the function $f(Z):=\sqrt{Z^{\top}AZ} = \|A^{1/2}Z\|_2$ has
Lipschitz constant $\mnorm{A^{1/2}}_2 = \sqrt{\mnorm{A}_2}$ with
respect to the Euclidean norm on $\R^p$. By the Gaussian isoperimetric
inequality~\cite{Massart}, we have for any $s>0$,
\begin{equation}
\label{gaussianIso}
\P\left[ f(Z) \leq \E[f(Z)] -s\right] \leq
\exp\left(\textstyle{\frac{-s^2}{2\mnorm{A}_2}}\right).
\end{equation}
From the Poincar\'e inequality for Gaussian measures~\cite{Beckner},
the variance of $f(Z)$ is bounded above as $\var[f(Z)] \leq
\mnorm{A}_2$.  Noting that $\E[f(Z)^2] = \tr(A)$, we see that the
expectation of $f(Z)$ is lower bounded as
\begin{equation*}
\E[f(Z)] \geq \sqrt{ \tr(A)-\mnorm{A}_2}.
\end{equation*}
Substituting this lower bound into the concentration
inequality~\eqref{gaussianIso} yields
\begin{equation*}
\P\left[ f(Z) \leq \sqrt{ \tr(A)-\mnorm{A}_2} -s\right] \leq
\exp\left(\textstyle{\frac{-s^2}{2\mnorm{A}_2}}\right).
\end{equation*}
Finally, letting $t\in
\left(0,\sqrt{\frac{\tr(A)}{\mnorm{A}_2}-1}\right)$, and choosing
$s^2= t^2\mnorm{A}_2$ yields the claim \eqref{upper}.\\

The Gaussian isoperimetric inequality also implies $\P\left[ f(Z)
  \geq \E f(Z) + s\right] \leq
\exp\left(\textstyle{\frac{-s^2}{2\mnorm{A}_2}}\right)$.  By Jensen's
inequality, we have
\begin{equation*}
\E[f(Z)] = \E\sqrt{Z^{\top} A Z} \leq \sqrt{\E[ Z^{\top} A Z] } =
\sqrt{\tr(A)},
\end{equation*}
from which we obtain $\P\left[ f(Z) \geq \sqrt{\tr(A)} + s\right] \leq
\exp\left(\textstyle{\frac{-s^2}{2\mnorm{A}_2}}\right)$, and setting
$s^2 = t^2\mnorm{A}_2$ for $t>0$ yields the claim \eqref{lower}.
\section{Proof of Proposition~\ref{PropKeyLimit}}
\label{AppPropKeyLimit}

The proof of Proposition~\ref{PropKeyLimit} is based on Lemmas~\ref{tracelemma}
and~\ref{frobeniuslemma}, which we state and prove below in Section~\ref{B2}. We then prove
Proposition~\ref{PropKeyLimit} in two parts, by first proving the
lower bound~\eqref{lowerProp}, and then the upper
bound~\eqref{upperProp} in sections \ref{lowerProof} and \ref{upperProof} respectively.
\subsection{Two auxiliary lemmas}\label{B2}
Note that the following two lemmas only deal with the randomness in the $k\times p$
matrix $P_k^{\top}$, and they can be stated without reference to the sample size $n$.

\begin{lemma}
\label{tracelemma}
Let $P_k^{\top}\in \R^{k\times p}$ have
i.i.d. $N(0,1)$ entries, where $k\leq p$. Assume there is a constant $a\in [0,1)$ 
such that $k/p\to a$ as $(k,p)\to\infty$. Then, there is a sequence of numbers $c_k\to
(1-\sqrt{a})^2$ such that
$$ \P\left[\frac{1}{p} \tr(P_k (P_k^{\top}\Sigma P_k)^{-1} P_k^{\top})
  \geq \frac{k}{\tr(\Sigma)} c_k\right]\to1 \text{ as } (k,p)\to\infty.$$
\end{lemma}
\emph{Proof.}  By the cyclic property of trace and Lemma
\ref{traceprod}, we have
\begin{align}
\label{secondstep}
\frac{1}{p} \tr\left(P_k(P_k^{\top}\Sigma P_k)^{-1}P_k^{\top}\right)&=
\frac{1}{p} \tr\left((P_k^{\top}\Sigma P_k)^{-1}P_k^{\top}
P_k\right)\\ & \geq \frac{1}{p} \tr\left((P_k^{\top}\Sigma
P_k)^{-1}\right) \lambda_{\min}\left(P_k^{\top} P_k\right).
\end{align}
For a general positive-definite matrix $A \in \R^{k\times k}$,
Jensen's inequality implies $\tr(A^{-1})\geq k^2/\tr(A)$. Combining
this with the lower-bound on $\lambda_{\min}(P_k^{\top}P_k)$ from
Lemma \eqref{wishart} leads to
\begin{equation}\label{r2}
\frac{1}{p} \tr\left(P_k(P_k^{\top}\Sigma P_k)^{-1}P_k^{\top}\right)
\geq \frac{k^2}{\tr(P_k^{\top}\Sigma P_k)}\cdot
\underbrace{[1-\sqrt{k/p}-t_2]^2}_{=:1-r_2(t_2)},
\end{equation}							      
with probability at least $1-\exp(-pt_2^2/2)$.\\ \\ We now obtain a high-probability upper bound on
$\tr(P_k^{\top}\Sigma P_k)$ that is of order $k \,\tr(\Sigma)$. First let $\Sigma= U\Lambda
U^{\top}$ be a spectral decomposition of $\Sigma$. Writing
$P_k^{\top}\Sigma P_k$ as $(P_k^{\top}U)\Lambda (U^{\top}P_k)$, and
recalling that the columns of $P_k$ are distributed as $N(0,I_{p\times
  p})$, we see that $P_k^{\top}\Sigma P_k$ is distributed as
$P_k^{\top}\Lambda P_k$. Hence, we may work under the assumption that
$\Sigma$ and $\Lambda$ are interchangeable. Let $0<\lambda_1\leq \cdots \leq \lambda_p$ be the eigenvalues of
$\Sigma$, with $\lambda_i=\Lambda_{ii}$, and let ${\bf{Z}} \in
\R^{(pk)\times 1}$ be a concatenated column vector of $k$ independent
and identically distributed $N(0,I_{p\times p})$ vectors. Likewise,
let ${\bf{\Lambda}}\in \R^{pk\times pk}$ be a diagonal matrix obtained
by arranging $k$ copies of $\Lambda$ along the diagonal, i.e.
\begin{align}
\bf{\Lambda} & \mydef \left ( \begin{array}{ccc} \Lambda & & \\ &
  \ddots & \\ & & \Lambda \end{array} \right).
\end{align}
By considering the diagonal entries of $P_k^{\top} \Lambda P_k$, it is
straightforward to verify that $\tr(P_k^{\top}\Lambda
P_k)\stackrel{d}{=} \bf{Z}^{\top} \bf{\Lambda} \bf{Z}.$
Applying Lemma \ref{quadform} to the quadratic form $\bf{Z}^{\top}
\bf{\Lambda} \bf{Z},$ and noting that $\mnorm{\Lambda}_F/\tr(\Lambda)$
and $\mnorm{\Lambda}_2/\tr(\Lambda)$ are at most 1, we have
\begin{equation}\begin{split}
 \tr(P_k^{\top}\Lambda P_k)&\leq \tr({\bf{\Lambda}}) +
 2\sqrt{t_3}\small{\mnorm{\bf{\Lambda}}_F}+2t_3\mnorm{{\bf{\Lambda}}}_2\\ &=
 k \tr(\Lambda) + 2\sqrt{t_3}
 \sqrt{k}\mnorm{\Lambda}_F+2t_3\mnorm{\Lambda}_2\\ &
 \footnotesize{\leq k \tr(\Lambda)\underbrace{\left(
     1+\frac{2\sqrt{t_3}}{\sqrt{k}}+\frac{2t_3}{k}\right)}_{=:1+r_3(t_3)},}\\
\end{split}
\end{equation}
with probability at least $1-\exp(-t_3)$, giving the desired upper bound on 
$\tr(P_k^{\top}\Lambda P_k)$. In order to combine the last
bound with \eqref{r2}, define the event
\begin{align*}
E_k & \mydef \left\{\frac{1}{p} \tr\left(P_k(P_k^{\top}\Sigma
P_k)^{-1}P_k^{\top}\right) \geq
\frac{k}{\tr(\Sigma)}\frac{1-r_2(t_2)}{1+r_3(t_3)}\right\},
\end{align*}
and then observe that $\P(E_k) \geq 1-\exp(-pt_2^2/2)-\exp(-t_3)$
by the union bound. Choosing $t_2=1/p^{1/4}$ and $t_3=\sqrt{k}$, we
ensure that $\P(E_k)\to 1$ as $(k,p) \to \infty$, and moreover,
that
\begin{equation*}
\frac{1- r_2(t_2)}{1+r_3(t_3)}\to (1-\sqrt{a})^2,\\
\end{equation*}
which completes the proof.  \qed
%
%
\begin{lemma}\label{frobeniuslemma}
Assume the conditions of Lemma~\ref{tracelemma}. Then for any $C >
\frac{(1+\sqrt{a})^2}{(1-\sqrt{a})^2}$, we have
\begin{align}
\P \left [\mnorm{P_k(P_k^{\top}\Sigma P_k)^{-1} P_k^{\top}}_F \leq
  \frac{ C \sqrt{k}}{\lambda_{\min}(\Sigma)}\right] & \to 1 \quad
\mbox{as $(k,p) \to \infty$.}
\end{align}
\end{lemma}
~\\ \emph{Proof.}  By the relation $\mnorm{A}_F^2 =
\tr(A^2)$ for symmetric matrices $A$, and the cyclic property
of trace,
\begin{align*}
\mnorm{P_k(P_k^{\top}\Sigma P_k)^{-1}P_k^{\top}}_F^2 &=
\tr\left(\left(P_k(P_k^{\top}\Sigma
P_k)^{-1}P_k^{\top}\right)^2\right) \; = \; \tr\left(
\left((P_k^{\top}\Sigma P_k)^{-1}P_k^{\top}P_k\right)^2\right).
\end{align*}
Letting $\rho(\cdot)$ denote the spectral radius of a matrix, we use
the fact that $|\tr(A)|\leq k \rho(A) \leq k \mnorm{A}_2$ for all real
$k\times k$ matrices $A$ (see \cite[p. 297]{HornJohnson}) to obtain
$$\mnorm{P_k(P_k^{\top}\Sigma P_k)^{-1}P_k^{\top}}_F^2 \leq k
\mnorm{\left((P_k^{\top}\Sigma P_k)^{-1}P_k^{\top}P_k\right)^2}_2.$$
Using the submultiplicative property of $\mnorm{\cdot}_2$ twice in
succession,
\begin{equation}\label{condnum}
\begin{split}
\mnorm{P_k(P_k^{\top}\Sigma P_k)^{-1}P_k^{\top}}_F^2 &\leq k
\mnorm{(P_k^{\top}\Sigma P_k)^{-1}P_k^{\top}P_k}_2^2\\ &\leq k
\mnorm{(P_k^{\top}\Sigma
  P_k)^{-1}}_2^2\cdot\mnorm{P_k^{\top}P_k}_2^2\\ &=
k\frac{1}{\lambda^2_{\min}(P_k^{\top} \Sigma P_k)}\cdot
\lambda_{\max}^2(P_k^{\top} P_k).
\end{split}
\end{equation}
Next, by Lemma~\ref{wishart}, we have the
bound
\begin{equation}\label{lmax} \lambda_{\max}(P_k^{\top} P_k) \leq
p\cdot\underbrace{[1+\sqrt{k/p}+t_4]^2}_{=:1+r_4(t_4)},
 \end{equation}
with probability at least $1-\exp(-pt_4^2/2)$.\\

By the variational characterization of eigenvalues, followed by
Lemma~\ref{wishart}, we have
\begin{equation}\label{lmin}
\begin{split}
\lambda_{\min}(P_k^{\top} \Sigma P_k)&=\inf_{\|x\|_2=1} \left(x^{\top}
P_k^{\top}\Sigma P_k x\right)\\ &\geq \inf_{\|y\|_2=1} \left(y^{\top}
\Sigma y\right) \inf_{\|x\|_2=1}\|P_k x\|_2^2\\ &=
\lambda_{\min}(\Sigma)\cdot \lambda_{\min}(P_k^{\top} P_k)\\ &\geq
\lambda_{\min}(\Sigma)\cdot p\cdot (1-r_2(t_2)),\\
\end{split}
\end{equation}
with probability at least $1-\exp(-pt_2^2/2),$ and $r_2(t_2)$ defined
as in line \eqref{r2}.\\ 

Substituting the bounds \eqref{lmax} and
\eqref{lmin} into line \eqref{condnum}, we obtain
\begin{equation}\label{lastevent}
\begin{split}
\mnorm{P_k(P_k^{\top}\Sigma P_k)^{-1}P_k^{\top}}_F^2  \leq \frac{k}{\lambda^2_{\min}(\Sigma)} \frac{ (1+r_4(t_4))^2}{(1-r_2(t_2))^2}.\\
\end{split}
\end{equation}
with probability at least $1-\exp(-pt_2^2/2)-\exp(-pt_4^2/2)$, where
we have used the union bound.\\

Setting $t_2=t_4=1/p^{1/4}$, the probability of the event
\eqref{lastevent} tends to 1 as $(k,p)\to\infty$. Furthermore,
$$\frac{ (1+r_4(t_4))^2}{(1-r_2(t_2))^2}\to
\frac{(1+\sqrt{a})^4}{(1-\sqrt{a})^4},$$ and so we may take $C$ in the
statement of the lemma to be any constant strictly greater than
$\frac{(1+\sqrt{a})^2}{(1-\sqrt{a})^2}$. \qed
\subsection{Proof of lower bound~\eqref{lowerProp} in Proposition~\ref{PropKeyLimit}}\label{lowerProof}
By the assumption on the distribution of $\delta$, we may write
$\delta/\|\delta\|_2 $ as $Z/\|Z\|_2$ where $Z\sim N(0,I_{p\times
  p})$. Furthermore, because $\|Z\|_2/\sqrt{p}\to 1$ almost surely as
$n\to\infty$, it is possible to replace $\delta/\|\delta\|_2$ with
$Z/\sqrt{p}$, and work under the assumption that
$\frac{\Delta_k^2}{\|\delta\|_2^2} = \frac{1}{p}
Z^{\top}P_k(P_k^{\top}\Sigma P_k)^{-1}P_k^{\top}Z$.  Noting that we
may take $Z$ to be independent of $P_k^{\top}$, the concentration
inequality for Gaussian quadratic forms in Lemma 2 
gives a lower bound on the
conditional probability
\begin{equation}
\label{firststep}
\P\left[ \textstyle \frac{\Delta_k^2}{\|\delta\|_2^2} \geq
  \frac{1}{p}\tr\left(P_k(P_k^{\top}\Sigma P_k)^{-1}P_k^{\top}\right)
  -\psi(t_1) \Big\bracevert P_k^{\top}\right] \geq 1-\exp(-t_1),
\end{equation}
where $\psi(t_1):= \frac{2\sqrt{t_1}}{p}\mnorm{P_k(P_k^{\top}\Sigma
  P_k)^{-1}P_k^{\top}}_F$ is a random error term, and $t_1$ is a
positive real number that may vary with $n$. Now that the randomness
from $\delta$ has been separated out in \eqref{firststep}, we study
the randomness from $P_k^{\top}$ by defining the event
\begin{equation}\label{event}
\textstyle
\mathcal{E}_n:= \left\{\frac{1}{p}\tr\left(P_k(P_k^{\top}\Sigma P_k)^{-1}P_k^{\top}\right) -\psi(t_1) \geq L_n\right\},
\end{equation}
where $L_n$ is a real number whose dependence on $n$ will be specified below. To see the main line of argument toward the statement of the proposition, we integrate the conditional probability in line~\eqref{firststep} with respect to $P_k^{\top}$, and obtain
\begin{equation}\label{Mainbound}
\P\left( \textstyle{\frac{\Delta_k^2}{\|\delta\|_2^2} \geq L_n} \right) \geq [1-\exp(-t_1)] \, \P(\mathcal{E}_n).
\end{equation}
The rest of the proof proceeds in two parts. First, we lower-bound $\tr\left(P_k(P_k^{\top}\Sigma P_k)^{-1}P_k^{\top}\right)\big/ p$ on an event $\mathcal{E}_n^{\prime}$ with $\P(\mathcal{E}_n^{\prime})\to 1$ as $n\to\infty$. Second, we upper-bound $\psi(t_1)$ on an event $\mathcal{E}_n^{\prime \prime}$ with $\P(\mathcal{E}_n^{\prime \prime})\to 1$. Then we choose $L_n$ so that $\mathcal{E}_n\supset\mathcal{E}_n^{\prime}\cap \mathcal{E}_n^{\prime \prime}$, and take $t_1\to\infty$ so that \eqref{Mainbound} implies $\P(\Delta_k^2\big/ \|\delta\|_2^2 \geq L_n) \to 1$ as $n\to\infty$.\\

For the first step of lower-bounding $\tr\left(P_k(P_k^{\top}\Sigma P_k)^{-1}P_k^{\top}\right)\big/p$, Lemma~\ref{tracelemma} asserts that there is a sequence of numbers $c_n\to (1-\sqrt{a})^2$ such that the event
\begin{equation}\label{traceevent}
\textstyle
\mathcal{E}_n^{\prime} :=\left\{ \frac{1}{p} \tr\left(P_k(P_k^{\top}\Sigma P_k)^{-1}P_k^{\top}\right) \geq \frac{k}{\tr(\Sigma)} c_n\right\}
\end{equation}
satisfies $\P(\mathcal{E}_n^{\prime})\to 1$ as $n\to\infty$.\\

Next, for the second step of upper-bounding the error $\psi(t_1)$, Lemma~\ref{frobeniuslemma} guarantees that for any constant $C$ strictly greater than $\frac{(1+\sqrt{a})^2}{(1-\sqrt{a})^2}$, the event 
\begin{equation}\label{frobeniusevent}
\textstyle
\mathcal{E}_n^{\prime \prime}:=\left\{  \frac{2}{p}\mnorm{P_k(P_k^{\top}\Sigma P_k)^{-1}P_k^{\top}}_F \leq \frac{ 2\, C \sqrt{k}}{p \lambda_{\min}(\Sigma)} \right\}
\end{equation}
satisfies $\P(\mathcal{E}_n^{\prime\prime})\to1$ as  $n\to\infty$.\\

 Now, with consideration to $\mathcal{E}_n^{\prime}$ and $\mathcal{E}_n^{\prime \prime}$, define the deterministic quantity
\begin{equation}\label{L}
\textstyle
L_n:= \frac{k}{\tr(\Sigma)} \left[ c_n - \sqrt{t_1} \frac{2\,C}{\sqrt{k}} \frac{\tr(\Sigma)}{p \lambda_{\min}(\Sigma)} \right],
\end{equation}
which ensures $\mathcal{E}_n\supset \mathcal{E}_n^{\prime} \cap
\mathcal{E}_n^{\prime \prime}$ for all choices of $t_1$. Consequently,
$\P(\mathcal{E}_n)\to 1$, and it remains to choose $t_1$ appropriately
so that the probability in line \eqref{Mainbound} tends to 1. If we
let $t_1 =\sqrt{k}\, \frac{p \, \lambda_{\min}}{\tr(\Sigma)}$, then
$t_1\to\infty$ by assumption ({\bf{A5}}), and the second term inside
the brackets in line \eqref{L} vanishes as $n\to\infty$. Altogether,
we have shown that
\begin{equation*}
L_n \frac{\tr(\Sigma)}{k} \to (1-\sqrt{a})^2, \quad \mbox{ and } \quad
\P\left( \textstyle{\frac{\Delta_k^2}{\|\delta\|_2^2} \geq L_n}
\right)\to 1.
\end{equation*}
It follows that $\P \left(\frac{\Delta_k^2}{\|\delta\|_2^2} \geq
  \frac{c\, k}{\tr(\Sigma)} \right) \to 1$ for any positive constant $c<(1-\sqrt{a})^2$, 
 which completes the proof of the
lower bound \eqref{lowerProp}.\qed\\
%
%

\subsection{Proof of upper bound~\eqref{upperProp} in Proposition~\ref{PropKeyLimit}}
\label{upperProof}
 As in the proof of the lower bound~\ref{lowerProp} in Appendix~\ref{lowerProof}, we may reduce to the case that
 $\frac{\Delta_k^2}{\|\delta\|_2^2} = \frac{1}{p}
 Z^{\top}P_k(P_k^{\top}\Sigma P_k)^{-1}P_k^{\top}Z$.  Conditioning on
 $P_k^{\top}$, Lemma \ref{quadform} gives a lower bound on the
 conditional probability
\begin{equation}\label{firststepUpper}
\P\left[ \textstyle \frac{\Delta_k^2}{\|\delta\|_2^2} \leq
  \frac{1}{p}\tr\left(P_k(P_k^{\top}\Sigma P_k)^{-1}P_k^{\top}\right)
  +\psi(s_1)+\phi(s_1) \Big\bracevert P_k^{\top}\right] \geq
1-\exp(-s_1),
\end{equation}
where $s_1$ is a positive real number that may vary with $n$, and we define
\begin{equation}
\psi(s_1) := \frac{2\sqrt{s_1}}{p}\mnorm{P_k(P_k^{\top}\Sigma
  P_k)^{-1}P_k^{\top}}_F, \quad \text{ } \quad \phi(s_1):=\frac{2
  s_1}{p}\mnorm{P_k(P_k^{\top}\Sigma P_k)^{-1}P_k^{\top}}_2.
\end{equation}
 Clearly,
$\phi(s_1)\leq \sqrt{s_1}\psi(s_1)$. Again, as in the proof of the
lower bound \eqref{lowerProp}, we let $U_n$ denote an upper bound
whose dependence on $n$ will be specified below, and we define an
event
\begin{equation}
\mathcal{D}_n:= \left\{ \frac{1}{p}\tr\left(P_k(P_k^{\top}\Sigma
P_k)^{-1}P_k^{\top}\right) +(1+\sqrt{s_1})\, \psi(s_1) \leq
U_n\right\},
\end{equation}
and integrate with respect to $P_k^{\top}$ to obtain
\begin{equation}\label{combine}
\P \left(\frac{\Delta_k^2}{\|\delta\|_2^2} \leq U_n \right) \geq
[1-\exp(-s_1)] \, \P(\mathcal{D}_n).
\end{equation}
Continuing along the parallel line of reasoning, we upper-bound
$\frac{1}{p}\tr\left(P_k(P_k^{\top}\Sigma P_k)^{-1}P_k^{\top}\right)$
on an event $\mathcal{D}_n^{\prime}$ (defined below) with
$\P(\mathcal{D}_n^{\prime})\to 1$, and re-use the upper bound of
$\psi(s_1)$ on the event $\mathcal{E}_n^{\prime \prime}$ (see line \eqref{frobeniusevent}), which was
shown to satisfy $\P(\mathcal{E}_n^{\prime \prime})\to
1$. Then, we choose $U_n$ so that $\mathcal{D}_n\supset
\mathcal{D}_n^{\prime} \cap \mathcal{E}_n^{\prime \prime}$, yielding $\P(\mathcal{D}_n)\to 1$.
Lastly, we take
$s_1\to \infty$ at an appropriate rate so that the probability in line
\eqref{combine} tends to 1.\\

To define the event
$\mathcal{D}_n^{\prime}$ for upper-bounding
$\frac{1}{p}\tr\left(P_k(P_k^{\top}\Sigma P_k)^{-1}P_k^{\top}\right)$,
note that for a symmetric matrix $A$ with rank $k$, Jensen's inequality implies
$\tr(A) \leq \sqrt{k \tr(A^2)}$, regardless of the size of
$A$. Considering $A = P_k(P_k^{\top}\Sigma P_k)^{-1}P_k^{\top}$, and
$\sqrt{\tr(A^2)} = \mnorm{A}_F$, we see that we may choose
$\mathcal{D}_n^{\prime}= \mathcal{E}_n^{\prime \prime}$ from line \eqref{frobeniusevent},
and on this set we have the inequality,
\begin{equation}
\frac{1}{p} \tr(P_k(P_k^{\top}\Sigma P_k)^{-1}P_k^{\top}) \leq
\sqrt{k}\, \frac{C\sqrt{k}}{p \, \lambda_{\min}(\Sigma)},
\end{equation}
with probability tending to 1 as $n\to\infty$, as long as $C$ is strictly
greater than $\frac{(1+\sqrt{a})^2}{(1-\sqrt{a})^2}$.
In order to guarantee the inclusion $\mathcal{D}_n\supset
\mathcal{D}_n^{\prime} \cap \mathcal{E}_n^{\prime
  \prime}$, we define
\begin{equation}\label{Udef}
U_n := \frac{C\, k}{p\, \lambda_{\min}(\Sigma)}\left[1 +
  (s_1+\sqrt{s_1})\,\frac{2}{\sqrt{k}}\right].
\end{equation}
Note that $k=\lfloor n/2\rfloor$ implies $k\to\infty$ as $n\to\infty$, so choosing $s_1=k^{1/4}$ ensures that $s_1\to\infty$ and the second term inside the brackets in line \eqref{Udef} vanishes. Combining lines \eqref{combine} and \eqref{Udef}, we have
\begin{equation*}
\P \left( \frac{\Delta_k^2}{\|\delta\|_2^2} \leq U_n\right) \to 1, \quad
\mbox{ and } \quad U_n \, \frac{p \, \lambda_{\min}(\Sigma)}{C\, k}\to 1.
\end{equation*}
It follows that $\P \left(\frac{\Delta_k^2}{\|\delta\|_2^2} \leq
\frac{C \, k}{p \, \lambda_{\min}(\Sigma)}\right) \to1 $ for any
constant $C$ strictly greater than
$\frac{(1+\sqrt{a})^2}{(1-\sqrt{a})^2}$, which completes the proof of
the upper bound \eqref{upperProp}.\qed
\bibliographystyle{unsrt}

\end{document}